\documentclass[journal,twoside,web]{ieeecolor}
\usepackage{generic}
\usepackage{cite}
\usepackage{amsmath,amssymb,amsfonts}
\usepackage{algorithmic}
\usepackage{graphicx}
\usepackage{textcomp}

\usepackage{eurosym}

\usepackage[T1]{fontenc} % Use 8-bit encoding that has 256 glyphs

\usepackage[utf8]{inputenc} % Required for including letters with accents

\usepackage{array,booktabs}
\usepackage{subfig}
\usepackage{makecell}
\usepackage{float}
\usepackage{pdfpages}
\usepackage{subfig} % Required for creating figures with multiple parts (subfigures)

\usepackage{amsmath,amssymb} % For including math equations, theorems, symbols, etc
\usepackage{mathtools}
\usepackage{amsfonts}
\usepackage{bbm}
\usepackage{multirow}
\usepackage{lastpage}
\usepackage{scrtime}
\usepackage{mathrsfs,dsfont} 
\usepackage{eurosym}
\usepackage{amscd}
\usepackage{caption}
\usepackage{varioref} % More descriptive referencing
\usepackage{tikz}
\usetikzlibrary {positioning}
\usepackage{pgfplots}
\pgfplotsset{ % Here we specify options for all figures in the document
	compat=newest, % Which version of pgfplots do we want to use?
	legend style =
	{font=\small \sffamily},
	label style = {font=\small\sffamily},
	every tick label/.append style={font=\small}}
\newcommand{\mc}{\mathcal}
\usepackage[most,skins,theorems]{tcolorbox}
\tcbset{highlight math style={enhanced,
		colframe=red,colback=white,arc=0pt,boxrule=1pt}}
\usetikzlibrary{arrows}

\def\1{\mathds{1}}

\usetikzlibrary[topaths]
\newcount\mycount

\DeclareMathOperator*{\rank}{rank}

\newtheorem{theorem}{Theorem}
\newtheorem{lemma}{Lemma}
\newtheorem{proposition}{Proposition}

\newtheorem{example}{Example}
\newtheorem{definition}{Definition}
\newtheorem{remark}{Remark}
\newtheorem{assumption}{Assumption}

\usepackage{url}
\newcommand{\ba}{\begin{aligned}}
\newcommand{\ea}{\end{aligned}}
\newcommand{\be}{\begin{equation}}
\newcommand{\ee}{\end{equation}}
\newcommand{\R}{\mathbb{R}}
\newcommand{\N}{\mathbb{N}}

\newcommand{\black}{\color{black}}
\newcommand{\tcb}{\textcolor{black}}
\begin{document}
	
	\title{Local identifiability of networks with nonlinear node dynamics}
	
	\author{
	Martina Vanelli,  ~\IEEEmembership{Member,~IEEE,}, and Julien M. Hendrickx
	\thanks{An earlier version of this paper {\black appeared at the 2025 European Control Conference (ECC) as \cite{Vanelli2025}.} This work was supported by F.R.S.-FNRS via the "KORNET" project, by the Concerted Research Action (ARC) via the
		“SIDDARTA” project % granted by the UCLouvain 
		and by FNRS via the “InterpoControl” Research Project.}
%	 \thanks{Martina Vanelli and Julien M. Hendrickx are with the Institute for Information and Communication Technologies, Electronics and Applied Mathematics (ICTEAM) UCLouvain, 1348 Ottignies-Louvain-la-Neuve, Belgium, email: \{martina.vanelli, julien.hendrickx\}@uclouvain.be. }
\thanks{{\black Martina Vanelli is with ENTEG and the J. C. Willems Center for Systems and Control, University of Groningen, Groningen, Netherlands, email:  {\tt \small m.vanelli@rug.nl}.} Julien M. Hendrickx {\black is} with the Institute for Information and Communication Technologies, Electronics and Applied Mathematics (ICTEAM) UCLouvain, 1348 Ottignies-Louvain-la-Neuve, Belgium, email:  {\tt \small julien.hendrickx@uclouvain.be}. }}
	\maketitle
	\thispagestyle{empty}
%	\maketitle
	
	\begin{abstract}
		We study the identifiability of nonlinear network systems with partial excitation and partial measurement when the network dynamics is linear on the edges and nonlinear on the nodes. We assume that the graph topology and the nonlinear functions at the node level are known, and we aim to identify the weight matrix of the graph. Our main result is 
		that, {\black for almost all  static analytic nonlinearities that cross the origin}, %{\black DAG}s are generically locally identifiable by exciting sources and measuring sinks 
		{\black directed graphs} are generically locally identifiable {\black if and only if at least one node is excited in every source component of the condensation graph and at least one node is measured in every sink component.} %every source component of the condensation graph contains at least one excited node and every sink component contains at least one measured node}.
%		provided that the nonlinearities are static, analytic and cross the origin. 
This holds even when all other nodes remain unexcited and unmeasured and stands in sharp contrast to most findings on network identifiability requiring measurement and/or excitation of each node. The result applies to \tcb{homogeneous} feed-forward \tcb{and recurrent} artificial neural networks \tcb{and generalizes previous literature by considering a broader class of activations and architectures.}  %in particular to feed-forward artificial neural networks with no offsets and generalize previous literature by considering a broader class of functions and topologies.
	\end{abstract}
	
	
	\section{Introduction}
	
Networks of interconnected dynamical systems are widespread across various domains \cite{Bullo2018}. Analyzing these systems and developing control strategies require understanding the interconnections, typically modeled as edges of a graph. However, identifying these systems from partial excitations and partial measurements is challenging because measured signals reflect combined dynamics \cite{Goncalves2008}.

There are various questions emerging in network identification, including learning the topology, identifying the edges, inferring global/local dynamics and analyzing spectral properties of the graph. In this work, we study network systems where the dynamics is linear on the edges and nonlinear on the nodes and our goal is to identify the weight matrix of the graph.
%\begin{figure}
%	\includegraphics[width=.5\textwidth]{fig1}
%	\caption{Example of dynamics where nodes 1 and 2 are excited		and node 4 is measured {\red (update figure and example to cyclic graphs)}}
%		\label{fig:ex1}
%\end{figure}
Specifically, we consider %weakly connected % acyclic graphs (DAGs)
{\black directed graphs (digraphs)}, denoted with $\mathcal{G}=(\mathcal{N},\mathcal{E},W)$, with set of nodes $\mathcal{N}$, set of directed links $\mathcal{E}$, and weight matrix $W$ whose entries are such that $w_{ij} \ne 0$ if and only if $(i,j)\in\mathcal{E}$. We assume that the output of a node $i$ at time ${\black t}$ is
\begin{equation}
	y_{i}^{{\black t}} = f_{i} \big(\sum_{j} w_{ij} y_{j}^{{\black t}-1} \big) + b_{i} u_{i}^{{\black t}-1} \label{eq:1}
\end{equation}
where $f_{i}$ is a given node-specific nonlinear  {\black activation function} and $u_{i}^{{\black t}-1}\in\mathbb{R}$ is an external excitation signal. The coefficient $b_{i}$ indicates whether the node is excited, namely, $b_{i}=1$ for excited nodes and $b_{i}=0$ otherwise. For example, \tcb{for the network in Fig. \ref{fig:ex1}, $b_2=1$ and $b_i=0$ for all $i\neq 2$. If the network is initially at rest (i.e., $y^{{\black t}}_i=0$ for all $i$ for all ${\black t}\leq0$), the output of node $3$ at time ${\black t}=4$ is $y_3^{{\black t}}=f_3(w_{34}f_4(w_{43}f_3(w_{32}u_2^{{\black t}-4}))+w_{32}(f_2(w_{21}f_1(w_{12}u_2^{{\black t}-4}))+u_2^{{\black t}-2}))$.} 
%$$\ba y_3^k=f_3(&w_{34}f_4(w_{43}f_3(w_{32}u_2^{k-4}))\\&+w_{32}(f_2(w_{21}f_1(w_{12}u_2^{k-4}))+u_2^{k-2}))\ea$$}%\tcb{consider the network in Fig. \ref{fig:ex1} where node $2$ is excited (i.e., $b_2=1$ and $b_i=0$ for $i\neq 2$) and node $3$ is measured. Then, if the network is initially at rest (i.e., $y^k_i=0$ for all $i$ and for all $k<0$), the output of node $3$ at time $k=5$ is $y_3^k=f_3(w_{34}f_4(w_{43}f_3(w_{32}f_2(u_2^{k-5})))+w_{32}f_2(w_{21}f_1(w_{12}u_2^{k-4})))$.}
\begin{figure}
	\centering
	\begin{tikzpicture}[scale=0.7]
		\foreach \x/\name in {(0,0)/1, (3,0)/2, (6,0)/3, (9,0)/4}\node[shape=circle,draw](\name) at \x {\small\name};
		
		\path [->,draw] (2) edge node[above] {\small $w_{32}$} (3);
		
		\foreach \a/\b/\w in {1/2/$w_{21}$, 3/4/$w_{43}$}\path [->,draw, bend left] (\a) edge node[above] {\small\w} (\b);
		
		\foreach \a/\b/\w in {2/1/$w_{12}$, 4/3/$w_{34}$}\path [->,draw, bend left] (\a) edge node[below] {\small\w} (\b);
		\foreach \x/\n/\name in { (2.5,1.5)/0/$u_2^{{\black t}}$, (6.5,1.5)/6/$y_3^{{\black t}}$}\node[](\n) at \x {\small \name};
		\foreach \x/\name in {(0,0.8)/$f_1$, (3.2,0.8)/$f_2$, (5.8,0.8)/$f_3$,(9,0.8)/$f_4$ }\node[](\name) at \x {\small\name};
		\foreach \a/\b in {0/2,3/6}\path [->,draw,dashed] (\a) edge node[above] {} (\b);
	\end{tikzpicture}
	\caption{\tcb{Example of dynamics where node $2$ is excited
			and node $3$ is measured.}}
	\label{fig:ex1}
	\vskip-0.5cm
\end{figure}
%	$$\begin{aligned}	y_3^k=f_3(&w_{34}f_4(w_{43}f_3(w_{32}f_2(u_2^{k-5})))\\&+w_{32}f_2(w_{21}f_1(w_{12}u_2^{k-4})))\,.	\end{aligned}$$}%For example, for the network in Fig. \ref{fig:ex1} where node 1 and 2 are excited (i.e., $b_{i}=1$ for $i=1,2$ and $b_{i}=0$ for the others), the output of node 4 at time $k$ is given by $y_{4}^{k}=f_{4}(w_{43}f_{3}(w_{31}u_{1}^{k-3}+w_{32}u_{2}^{k-3}) + w_{41}u_{1}^{k-2})$.

In our setting, the topology of the graph $\mathcal{E}$ and the nonlinear functions $\{f_{i}\}_{i\in\mathcal{N}}$ are known, while the weight matrix $W$ is unknown. We are then interested in determining the identifiability of $W$ with partial excitation and partial measurement, that is, we want to determine if, given the sets of excited and measured nodes, there exist two different weight matrices leading to the same input-output behaviors. In this latter case, we say that the network is not identifiable as different weight matrices cannot be distinguished by experiments in which we excite and measure only such nodes.

Network identifiability of linear systems (when $f_{i}(x)=x$ for all $i$) with partial excitation and partial measurement has been the subject of recent research \cite{Hendrickx2018,VanWaarde2019,Bazanella2019,Legat2020,Cheng2021,Cheng2023,Mapurunga2023,Legat2024}. While some graph theoretical conditions exist for full measurement scenarios, a general solution remains elusive. Furthermore, real-world systems are mostly nonlinear. Identifiability of nonlinear systems is more challenging and far less studied in the literature. Recent work  {\black studies the identifiability of nonlinear network dynamics  when only the topology is known and consider both additive \cite{Vizuete2024a} and non-additive  \cite{Vizuete2024b,Vizuete2025} forms of unknown nonlinearities, showing that excitation and/or measurement of every node is necessary and deriving conditions for identifiability with full excitation \cite{Vizuete2024a}, partial excitation/measurement \cite{Vizuete2024b} and full measurement \cite{Vizuete2025}. In contrast, we assume a particular example of non-additive nonlinear systems having the structure in \eqref{eq:1} and further assume that the nonlinearities at the node level are known, allowing to drastically reduce the number of excited/measured nodes.}  Related work has also examined the structural identifiability of unknown linear coupling of nonlinear dynamical systems \cite{Karabutov2025} and parameter identifiability of networks composed of nonlinear systems with nonlinear coupling functions \cite{Verdiere2024}. %{\red Add references of revier + emphasize differences with Antoine and Renato} % in nonlinear interconnected systems
\tcb{Beyond network identifiability, observation strategies based on parameter-identifiability analysis have been proposed for biological networks \cite{wang2017parameter}, while robust identification methods for stochastic nonlinear parameter-varying systems \cite{morato2021robust} and observer-based estimation techniques for nonlinear dynamical systems subject to sensor faults \cite{djordjevic2022sensor} have also been developed. These works address challenges related to partial observations, uncertainty, and parameter estimation, but do not study the graph-theoretic identifiability of network weights considered here.}
%\tcb{\cite{wang2017parameter,morato2021robust,djordjevic2022sensor}}

%In this work, we differ from prior work by focusing on a particular example of non-additive nonlinear systems having the structure in \eqref{eq:1}.  Such 
{ \black Systems having the structure in \eqref{eq:1}} arise in several settings where nonlinear node dynamics interact through linear interconnections. Examples include continuous thresholds models \cite{Zhong2019}, that are continuous generalizations of linear threshold models \cite{granovetter1978threshold}, and nonlinear consensus problems \cite{Bizyaeva2022,Gray2018}. In these two settings, the action (resp. the opinion) of an agent at time $k$ depends in a nonlinear manner on the weighted sum of the actions (resp. the opinions) of his neighbors at time $k-1$. Identifying the weights allows to quantify the strengths of the interconnections among the agents, as long with the sign when coordinating and anti-coordinating agents coexist \cite{Lekamalage2026Identifiability}. Networks systems of the form in \eqref{eq:1} can be found also in network games, production networks, and models of financial interactions (see \cite{Acemoglu2015Networks}).

Our main finding is that %{\black DAG}s (DAGs) 
{\black a digraph is} locally identifiable for almost all weights and almost all analytic functions that cross the origin if and only if %all sources are excited and all sinks are measured. 
{\black at least one node is excited in every source component and one node is measured in every sink component of its condensation graph}. To make this precise, we introduce a novel notion of genericity for analytic functions, formally defining when a property holds for almost all analytic functions in a given class. To the best of our knowledge, the notion of almost all analytic functions is novel, having been previously introduced only in \cite{Vanelli2025}, with a related notion derived in \cite{Vizuete2025}. The necessary and sufficient condition of exciting {\black  source and measuring sink components of the condensation graph}, derived for both homogeneous and heterogeneous nonlinearities, stands in sharp contrasts with findings in identifiability of network systems, where exciting and/or measuring every node is proved to be necessary in both linear \cite{Hendrickx2018} and nonlinear \cite{Vizuete2024b} systems. In our setting, the presence of \tcb{a \emph{known nonlinear activation}} enables the identification of paths where some nodes remain unexcited and unmeasured, while the linear dynamics on the edges provide sufficient structure for network identification.

Our result applies in particular to {\black  feed-forward and recurrent homogeneous artificial neural networks (ANNs), i.e., ANNs with no offsets \cite{Nielsen2015}. In this context, the objective is to determine if, given a function $F$ and a nonlinearity $f$, it is possible to determine the architecture, weights, and biases of all {\black ANNs} that realize $F$. Identifiability of ANNs has been extensively studied, starting with shallow networks with $\tanh$ and other activation functions \cite{Sussmann1992,Albertini1993}, and later extending to deep networks, first for $\tanh$ \cite{fefferman1993recovering,fefferman1994reconstructing} and more recently for broader classes of sigmoidal functions \cite{Vlacic2021,Vlacic2022}, ReLU \cite{BonaPellissier2023}, and monomials \cite{Finkel2024,Usevich2025}. These results typically rely on techniques tailored to the specific nonlinearity and architecture. In contrast, we provide a unified framework for generic local identifiability of homogeneous networks that applies to almost all activations in an infinite-dimensional class, containing $\tanh$, polynomials, and modern activations as Gaussian Error Linear Unit and Swish,	and to arbitrary digraph topologies, including  cycles. Another goal of this work is to connect ANN identifiability with identifiability and observability of network dynamics; this connection has also been explored in recent work on ReLU networks \cite{yang2025local,yang2026local}.} %{\color{red} If time and space allow, include tanh, Gelu and Swish. Also, check $F^t(...)$ in the static mapping.}

A preliminary version of this work was presented in \cite{Vanelli2025} where we proved generic local identifiability of fully-connected feed-forward networks by only exciting sources and measuring sinks in the class of analytic functions that cross the origin. In this work, we extend the previous results to all %DAGs
{\color{black} digraphs} and we derive results for heterogeneous nonlinearities. The remainder of the paper is structured as follows. {\black In Section \ref{sec:prob}, we present the problem setting, temporarily focusing on homogeneous nonlinearities, and %and we derive a static input-output mapping that connects dynamical networks and {\black ANNs} (in Section \ref{ss:io_map}). We present some fundamental examples and we 
introduce all the relevant definitions (including the novel definition of genericity in the class of functions in Section \ref{ss:gen}).}
Our main result is stated in Section \ref{sec:main} and proved in Section \ref{sec:proof}. In Section \ref{sec:het}, we derive an analogous result for heterogeneous nonlinearities. Finally, we discuss our main assumptions in Section \ref{sec:ass} and we derive the conclusions in Section \ref{sec:conc}.
	
	\begin{figure}
	\centering
	\begin{tikzpicture}[scale=0.6]
		\foreach \x/\name/\label in {(0,0)/$C_1$/1, (3,0)/$C_2$/2}\node[shape=circle,draw](\label) at \x {\footnotesize\name};
		
		\path [->,draw] (1) edge node[above] {} (2);
	\end{tikzpicture}
	\caption{\tcb{Condensation graph for Fig. \ref{fig:ex1}, where $C_1=\{1,2\}$ and $C_2=\{3,4\}$ are a source and a sink component, respectively.}}
	\label{fig:cond_graph}
\end{figure}
	\subsubsection*{Notation}\label{ss:not}
	We denote the out- and in-neighborhood of a node $i$ in $\mathcal{N}$ respectively by $\mathcal{N}_{i}^{+}=\{j\in\mathcal{N}|(j,i)\in\mathcal{E}\}$ and $\mathcal{N}_{i}^{-}=\{j\in\mathcal{N}|(i,j)\in\mathcal{E}\}$. \tcb{A} walk of length $l$ from node $i$ to node $j$ is a finite sequence of nodes $(i=i_{0},i_{1},...,i_{l}=j)$ with $\black (i_{k},i_{k-1})\in\mathcal{E}$ for all $k=1,...,l$. A path is a walk $P=(i=i_{0},...,{\black i_{l}}=j)$ with $i_{k}\ne i_{k^{\prime}}$ for all $k\ne k^{\prime}$, except possibly $i_{0}=i_l$. For two nodes $i$ and $j$ in $\mathcal{N}$, we denote with $\mathcal{P}^{i\rightarrow j}$ the set of paths from $i$ to $j$ and with $\text{dist}(i,j)=\min\{|P|, P\in\mathcal{P}^{i\rightarrow j}\}$ the distance between $i$ and $j$. %{\black The diameter of a digraph $\mathcal{G}$ is the maximum distance between any pair of nodes, defined as $\text{diam}(\mathcal{G}) = \max_{i,j \in \mathcal{N}} \text{dist}(i,j)$.} 
	{\black A digraph is strongly connected if there exists a directed path between any pair of distinct nodes. A strongly connected component (SCC) of $\mathcal{G}$ is a maximal strongly connected subgraph. Let $\mathcal{N}_c = \{C_1, \dots, C_q\}$ denote the set of node sets corresponding to the $q$ disjoint SCCs of $\mathcal{G}$, such that $\bigcup_{k=1}^q C_k = \mathcal{N}$. The \emph{condensation graph} of $\mathcal{G}$ is defined as the digraph $\mathcal{G}_c = (\mathcal{N}_c, \mathcal{E}_c)$, where the nodes represent the SCCs of $\mathcal{G}$. A directed edge $(C_i, C_j) \in \mathcal{E}_c$ exists if $i \neq j$ and there is at least one edge $(u, w) \in \mathcal{E}$ with $u \in C_i$ and $w \in C_j$.  A SCC $C_i$ is called a source component if it has no incoming edges in $\mathcal{G}_c$, and a sink component if it has no outgoing edges in $\mathcal{G}_c$. %Since $\mathcal{G}_c$ is a finite DAG, 
	By definition, the condensation graph %$\mathcal{G}_c$ contains no directed cycles, i.e., it
	 is a directed acyclic graph (DAG). Thus, it possesses at least one source and one sink component (see Fig. \ref{fig:cond_graph}).}
	Given a graph $\mathcal{G}$ and a subset of nodes $S\subseteq\mathcal{N}$, the subgraph of $\mathcal{G}$ induced by $S$ is given by $\mathcal{G}[S]=(S,\mathcal{F},W|_{S\times S})$, where $\mathcal{F}=\{(i,j)\in\mathcal{E}:i,j\in S\}$. Specifically, for a {\black digraph} $\mathcal{G}$ and two nodes $i, j$ in $\mathcal{N}$, we denote with $\mathcal{G}^{i\rightarrow j}=\mathcal{G}[S^{i\rightarrow j}]$ the subgraph induced by the set of nodes that belong to a path from $i$ to $j$, that is, $S^{i\rightarrow j}=\{k\in\mathcal{N} : \exists P\in\mathcal{P}^{i\rightarrow j}, k\in P\}$. {\black For $F, \tilde F:\mathbb{R}^{n}\rightarrow\mathbb{R}$, 
 $F\equiv \tilde F$ denotes $F(x)=\tilde F(x)$,  $\forall x\in \mathbb{R}^n$. For $f,g:\R\rightarrow \R$, $f=o(g)$ denotes $\lim_{x \to 0} \frac{f(x)}{g(x)} = 0$.}
	\section{Problem Formulation}\label{sec:prob}
	
	\subsection{Identifiability}\label{ss:id}
	The identifiability analysis will be based on the standard notion of identifiability in system identification \cite{lennart1999system}, that is, we aim to determine whether two different set of weights can lead to the same global behavior. This is equivalent of requiring that the weights could, in principle, be uniquely recovered from perfect knowledge of the network's global input-output behavior. This condition is necessary for successful system identification, because if the model structure is not identifiable, no amount of high-quality data can guarantee recovery of the true weight parameters. 
	
	We consider the model class in \eqref{eq:1}, {\black where each node is equipped with an activation function belonging to %{\black	with analytic activation functions that pass through the origin, i.e., belonging to %we focus on nonlinearities that are static, analytic and cross the origin, i.e., 
% we assume that $f$ belongs to the class of functions	
\begin{equation}
	\mathcal{F}_{Z} := \{f \in \mc C^{\omega }, f(0)=0\}\,, \label{eq:9}
		\end{equation} 
where  $ C^{\omega }$ denotes the space of analytic functions from $\R$ to $\R$. Examples of commonly used activation functions in $\mathcal{F}_Z$ are the hyperbolic tangent $f(x)=\tanh(x)$, the Gaussian Error Linear Unit (GELU) $f(x)=x\Phi(x)$, where $\Phi$ denotes the standard normal cumulative distribution function, the Swish activation $f(x)=x/(1+e^{-x})$, and polynomial activations with vanishing constant term.
	We initially consider homogeneous nonlinearities $f_{i}=f$ for all nodes $i$, while the case with heterogeneous nonlinearities is addressed in Section \ref{sec:het}.} %We assume that the nonlinearity $f$ and the graph topology $\mc E$ are known,} %
We assume to know the nonlinearity $f$ and the graph topology $\mathcal{E}$
	 and we want to determine the identifiability of $W$ when only a subset of the nodes is excited (partial excitation) and only a subset of the nodes is measured (partial measurement). 	 %{\red Comment}

	%%%%%% comment %%%%%%%%%%%%
	%This assumption can be restrictive in practice (since some commonly used activation functions, such as Sigmoid, Softplus and Radial Basis Functions, are nonzero at zero), but it is needed for our technical analysis. Moreover, this assumption is closely related to the assumption that the network contains no biases, which is also limiting. We discuss these considerations further in Section \ref{sec:ass}. {\red Emphasize that several works are even more limiting by assuming specific nonlinearities, e.g., hyperbolic tangent and monomials, expecially when considering deep networks.} %We study local identifiability when the class of functions considered is made of all analytic functions that cross the origin, i.e.,}
		%	{\black We assume that the activation function $f$ belongs to a class of scalar analytic functions, denoted with $C^{\omega}$, and we consider a class of functions $\mathcal{F}\subseteq C^{\omega}$.		Throughout the paper, we focus on nonlinearities that are static, analytic and cross the origin, i.e., we assume that the activation function $f$ belongs to the class of functions	\begin{equation}			\mathcal{F}_{Z} := \{f:\mathbb{R}\rightarrow\mathbb{R} \mid f \text{ analytic in } \mathbb{R}, f(0)=0\}\,. \label{eq:9}		\end{equation}}
%		so that the origin is an equilibrium point of the network. 
	
	Following \cite{Vizuete2024a}, we denote with $\mathcal{N}^{e}$ and $\mathcal{N}^{m}$ in $\mathcal{N}$ the set of excited and measured nodes, respectively. Accordingly, we set $b_{i}=1$ for $i \in \mathcal{N}^{e}$ and $b_{i}=0$ for $i \in \mathcal{N}^{m}$. {\black For a node $i\in\mathcal{N}^{m}$,} we let $\mathcal{N}_{i}^{e,p}\subseteq\mathcal{N}^{e}$ denote the set of excited nodes that have a path to node $i$. %Due to the absence of memory on the edges and the nodes, in the absence of cycles, the output of a node $i$ depends only a finite number of past values of the inputs. {\black However, f
	{\black	For general digraphs containing cycles, signals can propagate indefinitely and the output of a node $i$ can depend on an infinite number of past values of the inputs. %leading to an infinite history of past inputs. 
		To establish a general framework that renders the input-output mapping finite for cyclic graphs, we make the following standard assumption \cite{nelles2020nonlinear, campi2002finite, pillonetto2023full,khosravi2023kernel}.
		\begin{assumption}%[Initially at rest]
			\label{ass1}The network is initially at rest at time $t=0$, %that is, $u_{i}^{t} = 0$ for all $t\le 0$, 
			that is, $y_{i}^{{\black t}} = 0$ for all ${\black t} \le 0$ and all $i \in \mathcal{N}$.\end{assumption}

		%Under Assumption \ref{ass1}, even if cycles are present, the number of inputs that appear in the function associated with the measurement of a node at a time $k > 0$ is strictly finite.} 
		Under Assumption \ref{ass1}, even in the presence of cycles,  the measurement of a node at a time ${\black t} > 0$ depends on a finite number of past inputs. For every measured node $i \in \mc N^m$ and every ${\black t}>0$, we define $F^{{\black t}}_i$ as the nonlinear mapping that captures the dependence of the measurement $y_i^k$	on the past inputs of the excited nodes having a directed path to $i$.}
		Then, if we measure node $i$ at time ${\black t}$, we obtain
	\begin{equation}
		y_{i}^{{\black t}} = u_{i}^{{\black t}-1} + F_{i}^{\black {\black t}}(u_{1}^{{\black t}-2}, \dots, u_{1}^{{\black t}-d_{1}{\black -1}}, \dots, u_{n_i}^{{\black t}-2}, \dots, u_{n_i}^{{\black t}-d_{n_i}-1}) \label{eq:2}
	\end{equation}
	where $1,\dots,n_i \in \mathcal{N}_{i}^{e,p}$ and $d_{n_i}$ denotes the maximum delay of excited node $n_{i}$ affecting the measurements. %{\black  For DAGs, $d_{n_i}$ is topologically bounded even without Assumption \ref{ass1}. For cyclic graphs, %it is bounded by the measurement time $k$, that is, 
%		$d_{e}\leq k-1$ for all $e$ in $\mc N^e$ where $k$ is the time of the measurement $F^k_e$.}
{\black Note that, under Assumption \ref{ass1}, $d_{e}\leq {\black t}-1$ for all $e$ in $\mc N^e$. }
%The function $F^k_i$ is implicitly defined by \eqref{eq:2} and the inputs in $F_i^k$ correspond to the nodes that have a path to the measured node. 
	
{\color{black}	\begin{example}
	\label{ex:1}
		Consider the graph $\mathcal{G}$ in Figure \ref{fig:ex1} with $\mc N^e=\{2\}$ and let $f_i=f$ for all $i$. Then, under Assumption \ref{ass1}, if we measure the node $3$ at time $k=4$, we obtain %$F_3^3=y_3^k$ with $y_3^k=f_3(w_{34}f_4(w_{43}f_3(w_{32}f_2(u_2^{k-5})))+w_{32}f_2(w_{21}f_1(w_{12}u_2^{k-4})))$.
		\be\label{eq:ex1}
		\ba
		F_{3}^{{\black t}}(u^{{\black t}-2}_{2}, u^{{\black t}-4}_{2})=f(&w_{34}f(w_{43}f(w_{32}u_2^{{\black t}-4}))\\&+w_{32}(f(w_{21}f(w_{12}u_2^{{\black t}-4}))+u^{{\black t}-2}_{2}))\,.
		\ea\ee
Observe that the function $F_{3}^{{\black t}}$ depends on the inputs of the excited node 2 with different delays, with $d_2=3= {\black t}-1$.
		%\begin{align}
		%	F_{4}^{\black k}(u_{1}^{k-3}, u_{2}^{k-3}) &= f(w_{43}f(w_{31}u_{1}^{k-3}+w_{32}u_{2}^{k-3})+w_{41}u_{1}^{k-2}) \label{eq:3}
		%\end{align}
		%Observe that the function $F_{4}$ depends on the inputs of the excited nodes 1 and 2 that have a path to the node 4. %Node 3 has a path to 4 but has no inputs, as it is not excited.
	\end{example}}

	%{\black The set $F^k$ represents how the outputs of measured nodes $y_k$ depend on the excitation signals from $0$ to $k$.} 		
	{\black
			Under Assumption \ref{ass1}, the measurements must be taken at a  sufficiently large time instant $t$ %the excitation signals to propagate through all the network and 
			to ensure that all edges on paths from excited nodes to measured nodes are crossed, so that all weights influence the output at least once. %and all weights %belonging to paths between  an excited and a measured node influence the output at least once. 
			%In cyclic digraphs, this requirement may additionally depend on the cycle structure. %A sufficient condition is $k\geq L_\text{max}+D_\text{max}$ where $L_\text{max}$ denotes the length of the longest directed cycle and $D_\text{max}$ the maximum path length between any excitation/measurement node pair. Smaller values of $k$ may also be sufficient depending on the graph structure. 
			A sufficient condition is to set ${\black t}$ greater than the length of the longest trail (i.e., walk with no repeated edges)  between any excitation/measurement node pair. However, smaller values of ${\black t}$ may also be sufficient. %, that is, the maximum walk with no repeated edges.
			For instance, for the graph in Fig. \ref{fig:ex1}, the maximum trail between the source $2$ and the sink $3$ is $(2,1,2,3,4,3)$, which has length $5$ and yields ${\black t}=6$, but ${\black t}=4$ is sufficient to cross all the edges (see \eqref{eq:ex1} in Ex. \ref{ex:1}).} 	Given a set of measured nodes $\mathcal{N}^{m}$, the set of measured functions $F^{\black t}(\mathcal{N}^{m})$ associated with $\mathcal{N}^{m}$ {\black at time $t$} is:
			\begin{equation}
			F^{\black t}(\mathcal{N}^{m}) := \{F_{i}^{\black t}, i\in\mathcal{N}^{m}\}. \label{eq:4}
			\end{equation}%In particular, we will consider\be k_\text{min}:=\text{diam}(\mc G)\,,\ee where $\text{diag}(\mc G)$ denotes the \emph{diameter} of the graph, that is, the maximum shortest path distance between any pair of vertices in the graph.
			
	As previously explained and similarly to{\black \cite{Hendrickx2018, Vizuete2024a}}, identifiability can be characterized by determining if knowledge of $F^{\black k}$ uniquely determines $W$. We say that a weight matrix $W$ is consistent with the edge set $\mathcal{E}$ if $w_{ij}\ne0$ only if $(i,j) \in \mathcal{E}$. Throughout, we shall denote with $F^{\black k}(\mathcal{N}^{m})$ the set of measured functions generated by a weight matrix $W$ consistent with $\mathcal{E}$ and with $\tilde{F}^{\black k}(\mathcal{N}^{m})$ the set generated by another matrix $\tilde{W}$ consistent with $\mathcal{E}$.
	
	\begin{definition}\label{def:id}%[Identifiability of a network]
		Consider a graph $\mathcal{G}=(\mathcal{N},\mathcal{E},W)$ with sets of excited and measured nodes $\mathcal{N}^{e}$ and $\mathcal{N}^{m}$ respectively, and let $\mathcal{F}$ be a class of functions. A network $\mathcal{G}$ is \emph{identifiable} in $\mathcal{F}$ if, for any given $f \in \mathcal{F}$, {\black there exists ${\black t}>0$ such that } $F^{\black t}(\mathcal{N}^{m})=\tilde{F}^{\black t}(\mathcal{N}^{m})$ implies $W=\tilde{W}$.
	\end{definition}
	
	If a network is not identifiable, different weight matrices can generate exactly the same global behavior, and therefore, recovering the weights from several experiments is impossible. On the other hand, if the network is identifiable, the weights can be distinguished from all others. Then, if the functions in $F(\mathcal{N}^{m})$ can be well approximated after sufficiently long experiments, it becomes feasible to approximate the weights through excitation of nodes in $\mathcal{N}^{e}$ and measurement of nodes in $\mathcal{N}^{m}$. 
	We have the following Lemma that establishes {\black the minimal conditions for $\mathcal{N}^{e}$ and $\mathcal{N}^{m}$.} %{\black, which is based on the notion of condensation graph (see Section \ref{ss:not}). } 
	%{\black \begin{definition}[Condensation Graph] Let $\mathcal{G} = (\mathcal{V}, \mathcal{E})$ be a digraph. The condensation graph $\mathcal{C} = (\mathcal{V}_{\mathcal{C}}, \mathcal{E}_{\mathcal{C}})$ is defined as follows. Each node $S \in \mathcal{V}_{\mathcal{C}}$ is a Strongly Connected Component (SCC) of $\mathcal{G}$. A directed edge $(S_i, S_j)$ exists in $\mathcal{E}_{\mathcal{C}}$ if and only if there exists at least one edge in the original graph $\mathcal{G}$ from some node $u \in S_i$ to some node $v \in S_j$.\end{definition} By construction, $\mathcal{C}$ is always a {\black DAG} (DAG), even if $\mathcal{G}$ contains cycles.}
	\begin{lemma}\label{lemma:min_cond}
		Let $\mathcal{G}$ be a %{\black DAG} 
		{
			\black digraph and $\mc C(\mc G)$ be its condensation graph.} %Given $\mathcal{F}$ be a class of functions. Then, t
		Then, {\black under Assumption \ref{ass1},} the network $\mathcal{G}$ is identifiable in the class {\black of functions} $\mathcal{F}$ only if {\black at least one node in each source component of $\mc C(\mc G)$  is excited and at least one node in each sink component of $\mc C(\mc G)$ is measured}.
	\end{lemma}
	\begin{proof}
		{\black Suppose a source component $C_{i}$ in $\mc C(\mc G)$ contains no excited nodes. Since no edges enter $C_i$, and no node in $C_i$ is excited, all node outputs $y_j^{{\black t}}$, $j \in \mc C_i$, remain zero for all ${\black t}$. Consequently,} %By definition of a source in the condensation graph, there are no edges entering $C_{i}$ from other components.}  Then, {\black all $y$ are identically zero in $C_i$. Consequently, } 
	for every ${\black t}>0$, the measurements $F^{\black t}(\mathcal{N}^{m})$ do not depend on the weights {\black involving nodes in this component} %{\black within this component and} of its outgoing edges 
		and, therefore, these weights are not identifiable because any weight will lead to the same $\black F^{t}(\mc N^m)$. %{\black Similarly, suppose a sink component $S_{snk}$ in $\mc C(\mc G)$ contains no measured nodes. By definition of a sink, there are no outgoing edges from $S_{snk}$ to other components. Therefore, } the weights {\black within this component and} of its incoming edges do not affect the measurements and are, therefore, not identifiable. %Any signal entering $S_{snk}$ via its incoming edges is "lost" within the component and never reaches a measurement point. Therefore, the weights of the edges entering $S_{snk}$ do not affect the available measurements and are not identifiable.}
	%	If a source is not excited, the measurements $F(\mathcal{N}^{m})$ do not depend on the weights of its outgoing edges and, therefore, these weights are not identifiable because any weight will lead to the same $F$. 
	Similarly, if a sink {\black component} is not measured, the weights %of its incoming edges
	{\black involving nodes in this component} do not affect the measurements and are, therefore, not identifiable.
	\end{proof}
	
	According to Lemma \ref{lemma:min_cond}, the minimal {\black number} of excited and measured nodes %{\black  for DAGs} Correspond to the sets of sources and sinks, respectively.
	{\black coincides with the number of sources and sinks in the condensation graph, respectively.}
	{ \black For DAGs, this condition reduces to  the sets of sources and sinks, respectively.  In this work, we will show that, under suitable generic conditions on the functions in \eqref{eq:9} and on the weights, this condition is in fact sufficient for (local) identifiability.  }

	\subsection{Input-Output maps}\label{ss:io_map}
	{\black Let $\mathcal{G}$ be a {\black digraph} with set of excited nodes $\mathcal{N}^{e}$ and measured nodes $\mathcal{N}^{m}$  and let Assumption \ref{ass1} hold.  }
For the upcoming analysis, we wish to express {\black for every ${\black t}>0$} the input-output map as a static function
	\begin{equation}
		F^{\black t}: \mathbb{R}^{N} \rightarrow \mathbb{R}^{M} \label{eq:5}
	\end{equation}
	where $N{\black \leq } \sum_{e\in\mathcal{N}^{e}}d_{e}$ 
	 is the sum of delays that affect the measurements and $M=|\mathcal{N}^{m}|$ is the number of {\black measured nodes}. This mapping will also explicit the fundamental connection with {\black ANNs}. For a given time step ${\black t}$ and scalar inputs $u_{e}^{{\black t}}\in\mathbb{R}$, define for each {\black excited node} $e$ the delay-history vector
	\begin{equation*}
		x^{(e)} := \begin{bmatrix}u_{e}^{{\black t}-2}\\ \vdots\\ u_{e}^{{\black t}-d_{e}-1}\end{bmatrix} \in \mathbb{R}^{d_{e}}.
	\end{equation*}
	The global history vector is obtained by stacking these blocks in a vector $x := \text{col}(x^{(e)})_{e\in \mathcal{N}^{e}} \in \mathbb{R}^{N}$. Then, the full measurement map in \eqref{eq:5} {\black at time ${\black t}$} is given by
	\begin{equation*}
		F^{\black t}(x) := [F_{m}^{\black t}(x)]_{m\in\mathcal{N}^{m}}.
		\end{equation*}
			
			\addtocounter{example}{-1}
			\begin{example}[cont'd]
				{\black Consider the dynamics in Figure \ref{fig:ex1} with $f_i=f$ for all $i$ and let ${\black t}=4$. According to \eqref{eq:ex1}, we find $F^{{\black t}}:\mathbb{R}^{2}\rightarrow\mathbb{R}$ with  $x=x^{(2)}=[u_{1}^{{\black t}-2}, u_{1}^{{\black t}-4}]^\top$ and $ F^{{\black t}}(x)=f(w_{34}f(w_{43}f(w_{32}x_2))+w_{32}(f(w_{21}f(w_{12}x_2))+x_1))\,.$}
				%$\ba F^4(x)=f_3\big(&w_{34}f_4(w_{43}f_3(w_{32}x_2))\\&+w_{32}(f_2(w_{21}f_1(w_{12}x_2))+x_1)\big)\,.\ea $$}%$F(x) = f(w_{43}f(w_{31}x_{2}+w_{32}x_{4})+w_{41}x_{1})$.}
			\end{example}
			\begin{figure}%\begin{figure*}%[!b]
				\centering
				\begin{tikzpicture}[scale=0.7]
					\foreach \x/\name/\sub in {(0,0)/1/$3_0$, (3,0)/2/$4_1$, (3,2)/3/$2_1$, (6,0)/4/$3_2$, (6,2)/5/$1_2$,
						%(9,0)/6/$2_3$, 
						(9,2)/7/$2_3$, 
						(9,0)/8/$4_3$}\node[shape=circle,draw](\name) at \x {\small\sub};
					
					\foreach \a/\b/\w in {1/2/$w_{34}$,
						1/3/$w_{32}$, 
						2/4/$w_{43}$,
						3/5/$w_{21}$,
						4/7/$w_{32}$,
						4/8/$w_{34}$,
						5/7/$w_{12}$			}\path [<-,draw] (\a) edge node[above] {\small\w} (\b);
					
					\foreach \x/\n/\name in { (4.5,3.5)/9/$u_2^{{\black t}-2}$, %(11,0)/10/$u_2^{k-4}$, 
						(11, 2)/11/$u_2^{{\black t}-4}$}\node[](\n) at \x {\small \name};
					
					\foreach \x/\name in {(0,0.8)/$f_3$, (3,2.8)/$f_2$, (3,0.8)/$f_4$, (6,2.8)/$f_1$, (6,0.8)/$f_3$,
						(9,0.8)/$f_4$, (9,2.8)/$f_2$}\node[](\name) at \x {\small\name}; %(9,0.8)/$f_4$}
				
				\foreach \a/\b in {3/9,7/11}\path [<-,draw,dashed] (\a) edge node[above] {} (\b);
			\end{tikzpicture}
			\caption{\tcb{Unfolded graph for Figure \ref{fig:ex1}}}
			
			\label{fig:unfolded}
			\end{figure} 

			The measurement of a {\black node} $m$ depends on the input of each {\black excited node} with different delays. {\black To make the signal propagation explicit in time, we construct a time-unfolded DAG, called \emph{unfolded graph} and denoted by $H_m^{\black t}(\mathcal{G}) = (\mathcal{N}_H, \mathcal{E}_H, W_H)$. This graph  is obtained by creating a copy of each node at every time step and connecting these copies according to the structure of the original graph. %In this way, the original digraph is "unfolded" over time, so that every node has multiple time-shifted copies corresponding to different delays. 
				For illustration, unfolding the graph from Figure \ref{fig:ex1} yields the DAG shown in Figure \ref{fig:unfolded}. 	%	\begin{example}[continued] 	Let us apply Algorithm 1 in [14] to the graph in Figure \ref{fig:ex1} and the measured node $m=4$. Then, we obtain the unfolded digraph.	\end{example}
			The exact algorithm, adapted from Algorithm 1 in \cite{Vizuete2024a}, proceeds as follows. First, a copy of the measured node $m$ is created and denoted by $m_0$. This represents the {\black measured} node at time ${\black t}$. Then, for each previous time step {\black ${\black k}=1,\dots, {{\black t}}-1$,}	a copy $i_{{\black k}}$ is created for each original node $i$ in $\mc N$, each associated with the corresponding excitation signal $u_{i}^{{\black t}-{\black k}-1}$. %$x_{\pi(i_k)}=u_{i}^{{\black t}-{\black k}-1}$, with $\pi : N_H \to \{1, \dots, N\} \cup \{0\}$ each unfolded node $i \in N_H$ to the corresponding position in the input vector $x$. 
			Edges in the unfolded graph are then defined according to $\mc G$: for any original edge $(i, j)$, an edge is included from $j_{{\black k}+1}$ to $i_{{\black k}}$, ensuring that the unfolded graph captures exactly how delayed signals propagate through the network. Note that some node copies may end up isolated (no incoming or outgoing edges); these nodes are discarded, as they have no effect on identifiability.}
			
			{\black 	In the unfolded digraph $H^{{\black t}}_m(\mathcal{G})$, the weight associated with the edge $(i_{{\black k}}, j_{{\black k}+1})$ corresponds to the original weight $w_{ij}$. Hence, the weight matrix $W_{H}$ of the unfolded graph is populated by repetitions of the entries of $W$. %Let $\Omega \subseteq \mathbb{R}^{|\mathcal{E}|}$ be the parameter space of the original weights, and $\Omega_{H} \subseteq \mathbb{R}^{|\mathcal{E}_{H}|}$ be the parameter space of the unfolded graph. 
				We define the linear parameter mapping:
				\be \label{eq:phi}
				\varphi: \mathbb{R}^{|\mathcal{E}|} \rightarrow \mathbb{R}^{|\mathcal{E}_{H}|}\,,\quad W_{H}=\varphi(W)
				\ee
				such that the weight of the unfolded edge $(i_{{\black k}},j_{{\black k}+1})$ is assigned as $(W_{H})_{i_{{\black k}}, j_{{\black k}+1}} =w_{ij}$. %Let $x$ denote the stacked vector of all relevant finite inputs $u_{i}^{k-t-1}$. 
			}
			We now have a static input-output mapping of the form \eqref{eq:5} and a DAG associated to it. 	 This can be expressed in a recursive form by setting {\black $F_{v}^{\black t}(x)=x_{\pi(v)}$  for every source $v\in \mc N_H$} and letting
			\begin{equation}
\black				\begin{aligned}
				F_{v}^{\black t}(x) &= f\big(\sum_{v'\in \mc N_H}{\black w_{vv'}^H}F_{v'}^{\black {\black t}}(x)\big){+\black x_{\pi(v)}} \\&= {\black f\big(\sum_{v'\in \mc N_H}(\varphi(W))_{vv'} F_{v'}^{t}}(x)\big)+x_{\pi(v)}\label{eq:6}
				\end{aligned}
			\end{equation}
{\black			for every non-source $v \in \mathcal{N}_H$,  where $x_{\pi(v)}$ selects the input injected at node $v$ (or $0$ if none).} In this way, $\black F_{m_0}^{\black t}(x)$
			coincides with the measurement of {\black node} $m$ {\black at time $t$}.  {\black	To emphasize the weight sharing and its dependence on the original $W$, we sometimes write the input-output mapping as $F^{{\black t}}(x; \varphi(W))$. We then have that $\mc G$ is identifiable  if and only if  $\exists {\black t}>0$ such that for any $W, \tilde{W}$ consistent with $\mathcal{E}$: 
			 	\be\label{eq:id2}F^{{\black t}}(x; \varphi(W)) = F^{{\black t}}(x; \varphi(\tilde{W})) \quad \forall x\quad \Rightarrow \quad W = \tilde{W}\,.\ee
			 	This formulation provides a direct connection between the identifiability problem for dynamical networks and the identifiability of \emph{homogeneous} ANNs:
			 	the map $F^{{\black t}}(\cdot;\varphi(W))$ is the input-output map of a feedforward neural network on $H_m^{{\black t}}(\mathcal G)$. Our framework covers the following architectures.

			 		\emph{Feedforward neural networks (FNNs).}
			 		If $\mathcal G$ is a layered acyclic graph and every directed path from an excited node to a measured node has the same length $L$, then only the delay ${{\black t}}=L+1$ contributes to the output. In this case, the unfolded graph $H_m^{{\black t}}(\mathcal G)$ is isomorphic to the original graph, $\varphi$ reduces to the identity map, and the network coincides with a \emph{layered FNN} of arbitrary depth.
			 %\emph{Networks with skip connections}.
			 		If $\mathcal G$ is acyclic but contains directed paths of different lengths, then several delays contribute simultaneously to the output. The unfolded graph separates these contributions into different branches, yielding architectures analogous to \emph{residual} or \emph{skip-connected networks}. Imposing equality constraints among the corresponding unfolded inputs recovers the standard setting where the same signal is injected across different delays.
			 		
			 		\emph{Recurrent neural networks (RNNs).}
			 		If $\mathcal{G}$ contains directed cycles, the unfolding procedure coincides with the usual unrolling-in-time representation of an RNN. In particular, the  recurrent neural network model
			 		$$			 		h_{{{\black t}}}=\sigma_h (W_h u_{{\black t}}+U_h h_{{{\black t}}-1})\,, \quad y_{t}=\sigma_y(W_yh_{{\black t}})
			 		$$
			 		can be obtained by considering a graph with $|u_{{\black t}}|$ excited source nodes %receiving the excitation $\textbf{u}^{k-1}$ 
			 		connected to $|h_{{\black t}}|$ non-excited non-measured nodes with self loops and activation $\sigma_h$, and a final layer of $|y_{{\black t}}|$ measurements with activation $\sigma_y$.  In Fig. \ref{fig:rnn}, we show the result for a simple network with $W_h=w_{21}$, $U_h=w_{22}$, $W_y=w_{32}$, and $\sigma_h=\sigma_y=f$. Under Assumption $1$, the measurement of node $3$ at time ${\black t}=4$ gives
			 		$$
			 		y_3^k = f \Big( w_{32} f \big( w_{22} f(w_{21} u_1^{{\black t}-4}) + w_{21} u_1^{{\black t}-3} \big) \Big)\,.
			 		$$
			 		By defining a time-shifted input sequence to account for the edge delays, this mirrors the unrolled standard RNN.} 

			 	\begin{figure}
	\centering
	\begin{tikzpicture}[scale=0.6]
		\foreach \x/\name in {(0,0)/1, (3,0)/2, (6,0)/3}
		\node[shape=circle,draw](\name) at \x {\small \name};
		
		\path (2) edge [loop below] node {\small $w_{22}$} (2);
		
		\foreach \a/\b/\w in {1/2/$w_{21}$, 2/3/$w_{32}$}\path [->,draw] (\a) edge node[above] {\small\w} (\b);
		
		\foreach \x/\n/\name in { (-2,0.7)/0/$u_1^{{\black t}}$, (8,0.7)/4/$y_1^{{\black t}}$}\node[](\n) at \x {\small \name};

		\foreach \x/\name in {(0,0.8)/$f$,(6,0.8)/$f$}\node[](\name) at \x {\small\name};
		\foreach \x/\name in { (3,0.8)/$f$}\node[](\name) at \x {\small\name};
		
		\foreach \a/\b in {0/1,3/4}\path [->,draw,dashed] (\a) edge node[above] {} (\b);
	\end{tikzpicture}
	\caption{\black Simple recurrent neural network}
	\label{fig:rnn}
\end{figure}
			\subsection{Some fundamental examples and local identifiability}\label{ss:local}
			Let us start by providing a simple example where, depending on the choice of the nonlinearity $f$, the network may or may not be identifiable by only exciting sources and measuring sinks {\black of the condensation graph}.
				
			\begin{example}\label{ex:path}
					{\black Consider the digraph $\mc G=(\mc N, \mc E, W)$ in Figure \ref{fig:path} with $\mc N=\{1,2\}$ and $\mathcal{E}=\{(2,1),(1,2)\}$. Let $\mathcal{N}^{e}=\mathcal{N}^{m}=\{1\}$, so that we excite and measure node $1$, and let Assumption \ref{ass1} hold. If we denote by $x=u_{1}^{{{\black t}}-3}$ the input of node $1$, we obtain that the measurement of node $1$ at time ${{\black t}}=3$ is \be\label{eq:F1} F_{1}^{{\black t}}(x)=f(w_{12}f(w_{21}x))\,.\ee} 
				Let $f(x)=a_{1}x+a_{2}x^{2}$ for some $a_{1}, a_{2}$ in $\mathbb{R}$. Then,
				\begin{align}
					\tcb{F_1^{{\black t}}}(x) &= w_{\black 12}w_{21}a_{1}^{2}x + w_{\black 12}a_{1}a_{2}(1+a_{1}w_{\black 12})w_{21}^{2}x^{2} \notag \\
					&\quad + 2{\black w_ {12}^{2}}w_{21}^{3}a_{1}a_{2}^{2}x^{3} + {\black w_{12}^{2}}w_{21}^{4}{\black a_{2}^{3}}x^{4} \label{eq:7}
				\end{align}
				{\black F}or $a_{1}=0$ and $a_{2}\ne0$, we find {\black $F_1^t(x)= w_{12}^{2}w_{21}^{4} a_{2}^{3}x^{4}$ yielding} 
				$\tcb{F_1^{{\black t}}}=\tcb{\tilde F_1^{{\black t}}}$ if and only if ${\black w_{ 12}^{2}}w_{21}^{4}={\black\tilde{w}_{ 12}^{2}}\tilde{w}_{21}^{4}$. {\black I}f we set $\tilde{w}_{21}=\frac{1}{\gamma}w_{21}$ and $\tilde{w}_{\black 12}=\gamma^{2}w_{\black 12}$ for any $\gamma\notin\{0,1\}$, we obtain $\tcb{F_1^{{\black t}}}=\tcb{\tilde F_1^{{\black t}}}$ and $W\ne\tilde{W}$. {\black Hence}, the network is not identifiable for $f(x)=a_{2}x^{2}$ with $a_{2}\ne0$. We obtain the same result for $a_{1}\ne0$ and $a_{2}=0$. {\black Let} $a_{i}\ne0$ for $i \in \{1,2\}${\black. B}y equating $\tcb{F_1^{{\black t}}}=\tcb{\tilde F_1^{{\black t}}}$, we obtain $4$ conditions, that is,
				%\begin{align*}
				\be \label{eq:syst_eq}
					\tcb{F_1^{{\black t}}=\tilde F_1^{{\black t}}}\;\Leftrightarrow \; %\\&
					\begin{cases}
						w_{\black 12}w_{21} = \tilde{w}_{\black 12}\tilde{w}_{21} \\
						w_{\black 12}w_{21}^2(1+w_{\black 12}) = \tilde{w}_{\black 12}\tilde{w}_{21}^2(1+\tilde{w}_{\black 12}) \\
						{\black w_{ 12}^{2}}w_{21}^3 = {\black \tilde w_{ 12}^{2}} \tilde{w}_{21}^3 \\
						{\black w_{ 12}^{2}}w_{21}^4 = {\black \tilde w_{ 12}^{2}} \tilde{w}_{21}^4\,.
					\end{cases}
				\ee
				%\end{align*}
				{\black Substituting} the third equation in the fourth equation {\black yields} $w_{21}=\tilde{w}_{21}$ and therefore $w_{\black 12}=\tilde{w}_{\black 12}$. Then, the only solution is $W=\tilde{W}$ and the network is identifiable.
			\end{example}

		\begin{figure}
	\centering	
	\begin{tikzpicture}[scale=0.6]
		\foreach \x/\name in {(0,0)/1, (3,0)/2}\node[shape=circle,draw](\name) at \x {\small\name};
		\foreach \a/\b/\w in {1/2/$w_{21}$}\path [->,draw, bend left] (\a) edge node[above] {\small\w} (\b);
		
		\foreach \a/\b/\w in {2/1/$w_{12}$}\path [->,draw, bend left] (\a) edge node[below] {\small\w} (\b);
		\foreach \x/\n/\name in { (-2.5,1)/0/$x=u_1^{{{\black t}}-3}$, (-2.5,-1)/6/$y_1^{{\black t}}={\black F_1^{{\black t}}}(x)$}\node[](\n) at \x {\small \name};
		\foreach \x/\name in {(0,0.7)/$f$, (3,0.7)/$f$}\node[](\name) at \x {\small\name};
		\foreach \a/\b in {0/1,1/6}\path [->,draw,dashed] (\a) edge node[above] {} (\b);
	\end{tikzpicture}
	\caption{\tcb{Graph considered in Example \ref{ex:path}}}%\tcb{On the left, a recurrent neural network. On the right, a cyclic graph with $n = 2$ nodes where weights arerepresented by the values on the edges.}}
\label{fig:path}
\end{figure}

		%	\begin{figure}
		%		\includegraphics[width=0.5\textwidth]{fig4}
		%		\caption{Graph considered in Example 4}
		%		\label{fig:nn}
		%	\end{figure}
			
			In Example \ref{ex:path}, we observed that, when $f(x)=a_{1}x+a_{2}x^{2}$, $a_{i}\ne0$, for $i \in \{1,2\}$, {\black the dynamics in Figure \ref{fig:path}} %a path graph with 3 nodes is identifiable by only exciting the source and measuring the sink.
			is identifiable by {\black by exciting and measuring only node $1$. } 
			Despite being a very simple example, this observation is fundamental as it is in contrast with the findings in the literature of identifiability of network systems. Indeed, it was proved in \cite{Hendrickx2018} for linear systems and in \cite{Vizuete2024a} for {\black unknown nonlinear dynamics} %nonlinear systems 
				that a necessary condition for identifiability is to excite and/or measure every node of the graph. In our setting, %the presence of the {\black known}
				{\black we extend the linear setting by introducing a \emph{known nonlinear activation}}  on the nodes {\black which allows} to identify paths where some nodes are neither excited or measured, as the superposition principle does not further apply. %On the other hand, the linear dynamics on the edges provide enough structure to identify the network. 
				That said, identifiability may fail for certain nonlinear functions, e.g., for monomials (see $f(x)=a_{2}x^{2}$ in Example \ref{ex:path}). {\black To sum up}, identifiability can be achieved without exciting or measuring all nodes, but it depends on the specific nonlinear function at the nodes.
			
			\begin{example}\label{ex:ann}
				Let us consider a graph with $n=4$ nodes and edge set as in Figure \ref{fig:ann}. Assume that $\mathcal{N}^{e}=\{1\}$ and $\mathcal{N}^{m}=\{4\}$. In this case, the {\black measurement of node $4$ at time ${\black t}=3$ is} %measured function of node 4 is given by
				\begin{equation*}
					{\black F_{4}^t}(x) = f(w_{42}f(w_{21}x)+w_{43}f(w_{31}x))
				\end{equation*}
				where $x=u_{1}^{{\black t}-3}$ denotes the input of node 1. Observe that, in the special case when $w_{21}=w_{31}$, the output of node 4 is
				\begin{equation*}
					{\black F_{4}^t}(x) = f((w_{42}+w_{43})f(w_{21}x))
				\end{equation*}
				This implies that, for some particular weight matrices, we can at most identify $w_{42}+w_{43}$ and therefore the network is never identifiable. Let us now make a second fundamental remark: for every weight matrix $W$, if we set %$\tilde{w}_{21}=w_{31}$, $\tilde{w}_{31}=w_{21}$, $\tilde{w}_{43}=w_{42}$ and $\tilde{w}_{42}=w_{43}$
				{\black \begin{equation}\label{eq:node_per}
\tilde{w}_{21}=w_{31}\,,\;\tilde{w}_{31}=w_{21}\,,\;\tilde{w}_{43}=w_{42}\,,\;\tilde{w}_{42}=w_{43}\,,\;
					\end{equation}
					we find} $W\ne\tilde{W}$ and ${\black F_{4}^{{\black t}}}={\black \tilde F_{4}^{{\black t}}}$. In words, we can at most identify edges up to node permutations. Therefore, in the general setting, the network is not identifiable.
			\end{example}
						\begin{figure}
				\centering
				\begin{tikzpicture}[scale=0.6]
					\foreach \x/\name in {(0,0)/1, (3,-1)/2, (3,1)/3, (6,0)/4}\node[shape=circle,draw](\name) at \x {\small\name};
					\foreach \a/\b/\w in {1/3/$w_{31}$,1/2/$w_{21}$,2/4/$w_{42}$,3/4/$w_{43}$}\path [->,draw] (\a) edge node[above] {\small\w} (\b);
					\foreach \x/\n/\name in { (-2.5,0.7)/0/$x=u_1^{{\black t}-3}$, (8,1)/6/$y^k_4={\black F_4^{\black t}}(x)$}\node[](\n) at \x {\small \name};
					\foreach \x/\name in {(0,0.7)/$f$,(3,-0.3)/$f$,(3,1.7)/$f$, (6,0.7)/$f$}\node[](\name) at \x {\small\name};
					\foreach \a/\b in {0/1,4/6}\path [->,draw,dashed] (\a) edge node[above] {} (\b);
				\end{tikzpicture}
				\caption{Graph considered in Example \ref{ex:ann}}
				\label{fig:ann}
			\end{figure}
			
			According to Example \ref{ex:ann}, the identifiability of the graph in Figure \ref{fig:ann} cannot be guaranteed unless we introduce more hypotheses and/or some relaxations. This problem is well-known in the context of identifiability of 
			{\black ANNs}\cite{Sussmann1992, Albertini1993}. In their setting, the weight matrix is assumed to satisfy some generic assumptions (e.g., weights within a layer must differ) to avoid anomalous behaviors (see irreducibility in \cite{Albertini1993}). Furthermore, identifiability is studied up to an equivalence class, accounting for transformations such as neuron permutations and sign flips {\black (for odd activations)},  {\black as well as rescaling transformations in ReLU \cite{BonaPellissier2023} and monomials \cite{Finkel2024,Usevich2025}.} %{\black  or rescaling (for ReLU \cite{BonaPellissier2023} and polynomial NNs \cite{Finkel2024,Usevich2025})} %to avoid the anomalous behavior observed for monomials in Example \ref{ex:path}).} The goal is then to fully characterize this class, ensuring no other transformations produce the same output.
			%In this work, instead, we consider the \tcb{notion} of local identifiability, which corresponds to identifiability provided that $\tilde{W}$ is sufficiently close to $W$, i.e., {\black $\|\tilde W-W\|<\epsilon$ for some given $\epsilon>0$.}
			{\black In their setting, the objective is to characterize the full equivalence class of parameters that yield the same input–output map. In this work, we instead adopt the notion of local identifiability, requiring uniqueness only within a sufficiently small neighborhood of the true parameters. Specifically, we consider perturbations $\tilde{W}$ satisfying $\|\tilde{W}-W\|<\epsilon$ for some $\epsilon>0$.}
%			{\black We remark that sign flips and node permutations are not local transformation, and therefore the network can be locally identifiable even if these transfromation give the same output. However, rescaling is a local trasnformation and therefore a network is not locally identifiable if rescaling with same output are possible. }			
			\begin{definition}\label{def:loc_id}%[Local identifiability]
				The network $\mathcal{G}=(\mathcal{N},\mathcal{E},W)$ is \emph{locally identifiable} in a class $\mathcal{F}$ with excitations $\mathcal{N}^{e}$ and measurement $\mathcal{N}^{m}$ if, for any given $f$ in $\mathcal{F}$, {\black there exist ${\black t}>0$ and $\epsilon >0$ such that $F^{{\black t}}(\mathcal{N}^{m}) = \tilde{F}^{{\black t}}(\mathcal{N}^{m})$ implies $W = \tilde{W}$ for any $\tilde{W}$ consistent with the graph satisfying $\|W-\tilde{W}\|<\epsilon$.} %there exist  $k>0$ and $\epsilon>0$ such that for any $\tilde{W}$ consistent with the graph satisfying $\|W-\tilde{W}\|<\epsilon$ there holds {\black $F^k(\mathcal{N}^{m}) = \tilde{F}^k(\mathcal{N}^{m})$ implies $W = \tilde{W}$.}
			\end{definition}
			
			{\black With this definition, the permutation ambiguity in Example \ref{ex:ann} does not necessarily prevent local identifiability. Indeed, if $w_{21}\neq w_{31}$ and $w_{42}\neq w_{43}$, the permuted parameter $\tilde W$ in \eqref{eq:node_per} is separated from $W$ by a positive distance and can therefore be excluded by choosing $\epsilon>0$ sufficiently small. 
				
	We remark that discrete and continuous symmetries have different implications for local identifiability. Discrete transformations, such as node permutations and sign flips, move $W$ to isolated equivalent parameters and therefore do not affect local identifiability. In contrast, continuous transformations, such as the rescaling
$\tilde w_{21}=\frac{1}{\gamma}w_{21}$, 				$\tilde w_{12}=\gamma^2w_{12}$ occuring in Example \ref{ex:path} for $f(x)=x^2$, produce equivalent parameters arbitrarily close to $W$ and hence prevent local identifiability. Finally, local identifiability may itself fail at particular parameter values (e.g.,  $w_{21}=w_{31}$ in Example \ref{ex:ann}). This motivates the introduction of \emph{generic identifiability}, as discussed below.}
			\subsection{Genericity}\label{ss:gen}
			The notion of generic identifiability has been already introduced for transfer functions in several works {\black (e.g., \cite{Hendrickx2018})}. A property holds generically, or for almost all variables, if it holds for all variables except a zero measure set in $\mathbb{R}^{N}$ where $N$ is the number of parameters of each variable. Since the edge set is known, the number of parameters of the weight matrix $W$ is $N=|\mathcal{E}|$, that is, we shall study identifiability of $W$ for almost all choices of the elements $w_{ij}$ that are not known to be zero. More precisely, we have the following definition.
			
			\begin{definition}\label{def:genW}%[\black Genericity in the weights]
				Given a graph $\mathcal{G}$ with sets $\mathcal{N}^{e}$ and $\mathcal{N}^{m}$ of excited and measured nodes and a fixed nonlinearity $f$, we say that the network is $W$-generically (locally) identifiable if it is (locally) identifiable for all $W$ consistent with $\mathcal{E}$ except possibly those lying on a zero measure set in $\mathbb{R}^{|\mathcal{E}|}$.
			\end{definition}
			
			In the remainder of the paper, we will consider a second type of genericity: genericity in the class of functions $\mathcal{F}$. {\black As observed in Example \ref{ex:path}, specific choices of nonlinearities, such as monomials, may lead to a failure of identifiability. We therefore investigate local identifiability from a generic perspective and prove that it holds for almost all functions in $\mathcal{F}_{\mathcal{Z}}$. In this sense, the failures exhibited in Example \ref{ex:path} are exceptional rather than typical.} %By establishing identifiability for almost all functions in $\mathcal{F}_{\mathcal{Z}}$, we avoid such pathological cases and obtain results that hold generically.} 
			We recall that genericity typically refers to all variables except possibly some lying on a lower dimensional space/zero measure set. In infinite dimension, these notions present some challenges. Here, we restrict the attention to the set of analytic functions, %, denoted with $C^{\omega}$, and we consider a class of functions $\mathcal{F}\subseteq C^{\omega}$. For analytic functions, indeed,
			{\black for which} we can make use of the series of the Maclaurin coefficients%{\black %(see Section \ref{ss:not})} 
			\footnote{Taylor coefficients in $x = 0$.} 
			to describe the function. We then define genericity by applying the zero measure set/lower dimensional subspace to a subset of Maclaurin coefficients. 
			\begin{definition}\label{def:genf}%[\black Genericity in the functions]
				Given a class of functions $\mathcal{F}\subseteq C^{\omega}$ a property holds $f$-generically in $\mathcal{F}$, or for almost all functions $f$ in $\mathcal{F}$, if  $\exists M>0$ such that the property holds for all functions $f$ in $\mathcal{F}$ except possibly a subset of functions whose first $M$ Maclaurin coefficients all lie on a zero measure set in $\mathbb{R}^{M}$.
			\end{definition}
			\addtocounter{example}{-2}
			\begin{example}[cont'd]%\label{ex:gen_id}
				{\black Consider %a path graph with 3 nodes (see
				the dynamics in Figure \ref{fig:path} under Assumption \ref{ass1}, and let $f \in \mathcal F_Z$} have Maclaurin coefficients $a=\{a_{k}\}_{k\in\mathbb{N}}$ {\black with $a_0=0$, that is,  there exists $r>0$ such that \be\label{eq:f_an}f(x)=\sum _{k \in \N_+}a_k x^k\,,%\,,\quad \forall x \in B_r(0): = \{x \in \R : \|x \| < r\}\,.
				\ee for all $x$ in $B_r(0): = \{x \in \R : \|x \| < r\}$.} %{\black Since $f$ is in $\mc F_{Z}$}, we let $a_{0}=0$. 
				Since $f$ is analytic, the measurement {\black $F_1^{{\black t}}$ in \eqref{eq:F1} is also analytic and there exists $\rho>0$, $\rho \leq r$,\footnote{\black If $f$ is analytic on $B_{r_1}(0)$ and $g$ is analytic on $B_{r_2}(0)$ with $g(0) = 0$, the composition $f \circ g$ is well-defined and analytic on $B_\rho(0)$, where $\rho = \sup \{ r \in (0, r_2) : g(B_r(0)) \subseteq B_{r_1}(0) \}$.}  such that for all $x \in B_\rho(0)$
				\begin{equation}
					F_{1}^{{\black t}}(x) = a_{1}^{2}w_{21}w_{12}x+a_{2}a_{1}w_{21}^{2}w_{12}(1+a_{1}w_{12})x^{2}+o(x^{2})\,, \label{eq:8}
				\end{equation}
				where the expansion follows by substituting \eqref{eq:f_an} in \eqref{eq:F1}. %By equating the first two coefficients of $F_{1}^k$ we obtaing the same system as in \eqref{eq:syst_eq}, which gives that} 
				 If $F_1^{{\black t}}\equiv \tilde F_1^{{\black t}}$, the first two Maclaurin coefficients of $F_1^t$ and $\tilde F^t$ must be equal, yielding $w_{12}=\tilde w_{12}$ and $w_{21}=\tilde w_{21}$ when $a_1a_2\neq 0$ as in \eqref{eq:syst_eq}. Hence,} the network is identifiable for all $f$ whose first two Maclaurin coefficients satisfy $a_{1}\ne0$ and $a_{2}\ne0$. {\black Since the set} $E=\{(a_{1},a_{2})\in\mathbb{R}^{2}:a_{1}a_{2}=0\}$ is a zero measure set of $\mathbb{R}^{M}$ with $M=2${\black, the} network is identifiable for almost all analytic functions satisfying $f(0)=0$. 
				 {\black Note that the condition $a_1a_2\neq 0$ is satisfied by any polynomial containing both linear and quadratic terms, by $\mathrm{GeLU}$ ($a_1=\frac{1}{2}$, $a_2=\frac{1}{\sqrt{2\pi}}$), and by $\mathrm{Swish}$ ($a_1=\frac{1}{2}$, $a_2=\frac{1}{4}$), but not by $\tanh(x)$, whose first two Maclaurin coefficients are $a_1=1$ and $a_2=0$\footnote{ \black The corresponding Maclaurin expansions are $\tanh(x)=x-\frac{x^3}{3}+o(x^5)$, $\mathrm{GeLU}(x)=\frac{1}{2}x+\frac{1}{\sqrt{2\pi}}x^2+o(x^4)$, and $\mathrm{Swish}(x)=\frac{1}{2}x+\frac{1}{4}x^2+o(x^3)$.}. Nevertheless, identifiability can be established for $\tanh(x)$ by expanding \eqref{eq:8} to the third coefficient and solving the corresponding system. Thus, the condition $a_1a_2\neq0$ is sufficient, but not necessary, for identifiability.} %{\black Note that identifiability is in particular satisfied by $\tanh(x)$ ($a_1=$ and $a_2=$), GeLU (), Swish and any polynomial with $a_1a_2\neq 0$.}
			\end{example}
		
			{	\black   
				\begin{remark}
					Observe that, in Example \ref{ex:path}, the first two Maclaurin coefficients of $F_{1}^{{\black t}}$ in \eqref{eq:8} depend \emph{only} on the first two Maclaurin coefficients  of $f$, namely, $a_{1}$ and $a_{2}$. More generally, we have that, for $f(0)=0$, the Maclaurin coefficients of the measurement $F^{{\black t}}$ depend only  a \emph{finite} number of Maclaurin coefficients of $f$. This property no longer holds when $f(0)\neq 0$, which motivates our assumption on the class $\mc F_{\mc Z}$. This point will be discussed in detail in the next sections, in particular in the proof of Theorem \ref{th:1} (see, e.g., Remark \ref{rem:f0} in Section \ref{ss:gen_W}) and in Section \ref{sec:ass}.  
					\end{remark}}
				%	{	\black We emphasize that the assumption that $f(0)=0$ guarantees that every Maclaurin coefficient of the measurement $F$ depends on a finite number of Maclaurin coefficients of the nonlinearity $f$. For instance, observe that in Example 5 the first two coefficients of $F_{3}$ in \eqref{eq:8} depend only on the coefficients $a_{1}$ and $a_{2}$. This is no longer true when $a_{0}\ne0$ (see Remark 2 for further details). }
		
			\section{Main Result}\label{sec:main}
	%		We study local identifiability when the class of functions considered is made of all analytic functions that cross the origin, i.e.,
	%		\begin{equation}
	%			\mathcal{F}_{Z} := \{f:\mathbb{R}\rightarrow\mathbb{R} \mid f \text{ analytic in } \mathbb{R}, f(0)=0\}. \label{eq:9}
	%		\end{equation}
	%		We assume that $f(0)=0$ since it guarantees that every Maclaurin coefficient of the measurement $F$ depends on a finite number of Maclaurin coefficients of the nonlinearity $f$. For instance, observe that in Example 5 the first two coefficients of $F_{3}$ in \eqref{eq:8} depend only on the coefficients $a_{1}$ and $a_{2}$. This is no longer true when $a_{0}\ne0$ (see Remark 2 for further details). This assumption can be restrictive in practice (since some commonly used activation functions, such as Sigmoid, Softplus and Radial Basis Functions, are nonzero at zero), but it is needed for our technical analysis. Moreover, this assumption is closely related to the assumption that the network contains no biases, which is also limiting. We discuss these considerations further in Section \ref{sec:ass}.		
			Our main result considers both forms of genericity introduced in Section \ref{ss:local}, according to the following definition.
			
			\begin{definition}\label{def:gen_Wf}%[\black Genericity in weights and functions]
				Given a graph $\mathcal{G}$ with sets $\mathcal{N}^{e}$ and $\mathcal{N}^{m}$ of excited and measured nodes and a class of functions $\mathcal{F}\subset C^{\omega}$, we say that the network is $f$-$W$-generically (locally) identifiable in the class $\mathcal{F}$ if it is $f$-generically $W$-generically (locally) identifiable in $\mathcal{F}$, that is, if it is $W$-generically (locally) identifiable for almost all functions in $\mathcal{F}$.
			\end{definition}
			
			We now state our main result, proved in Section \ref{sec:proof}.
			
			\begin{theorem}\label{th:1}
				Let $\mathcal{G}$ be a {\black digraph, with condensation graph $\mc C(\mc G)$,} and let $\mathcal{F}_{Z}$ be as in \eqref{eq:9}. 
					Then, {\black under Assumption \ref{ass1},} the network $\mathcal{G}$ is $f$-$W$-generically locally identifiable in the class $\mathcal{F}_{Z}$ if and only if  {\black at least one node in each source component of $\mc C(\mc G)$  is excited and at least one node in each sink component of $\mc C(\mc G)$ is measured}.
				%Let $\mathcal{G}$ be a {\black DAG} and let $\mathcal{F}_{Z}$ be as in \eqref{eq:9}. Then, the network $\mathcal{G}$ is $f$-$W$-generically locally identifiable in the class $\mathcal{F}_{Z}$ if and only if all sources are excited and all sinks are measured.
			\end{theorem}
			
			According to Theorem \ref{th:1}, %{\black for analytic functions that are zero in zero}, {\black digraphs} are $f$-$W$-generically locally identifiable by only exciting 
			{\black any digraph %with arbitrary topology and depth
				 is locally identifiable for almost all functions in $\mc F_Z$ and almost all weights consistent with the topology by only exciting one node in every source component and measuring one node in every sink component of its condensation graph. By Lemma \ref{lemma:1s1s}, this condition is also necessary, making it the weakest possible graph-theoretic requirement. This sharply contrasts with previous identifiability results where measuring/exciting each node is proved to be
				 a necessary condition  for both linear\cite{Hendrickx2018} and nonlinear network systems \cite{Vizuete2024a}, showing that both the presence and the knowledge of the nonlinear node activation fundamentally changes the identifiability problem.} %The result therefore shows that known nonlinear node dynamics fundamentally increase the amount of information that can be extracted from partial excitation and measurement.} 
			{\black Notably, Theorem \ref{th:1} applies to arbitrary  topologies, including graphs containing directed cycles, with arbitrary size and to almost all activation functions in the infinite-dimensional class $\mc F_Z$. Unlike previous identifiability results in ANNs, whose proofs rely on specific activation functions (e.g., sigmoidal nonlinearities \cite{Fefferman1994, Vlacic2021,Vlacic2022}, ReLU \cite{BonaPellissier2023}, monomials \cite{Finkel2024,Usevich2025}) or restricted architectures such as FNNs, our result provides a unified approach for generic local identifiability in a function class. The theorem relies on the assumptions that the activation functions belong to $\mc F_Z$ and that no bias terms are present. The mathematical role and practical implications of these assumptions are discussed in Section \ref{sec:ass}. }

\begin{figure}
			\centering	\includegraphics[width=0.45\textwidth]{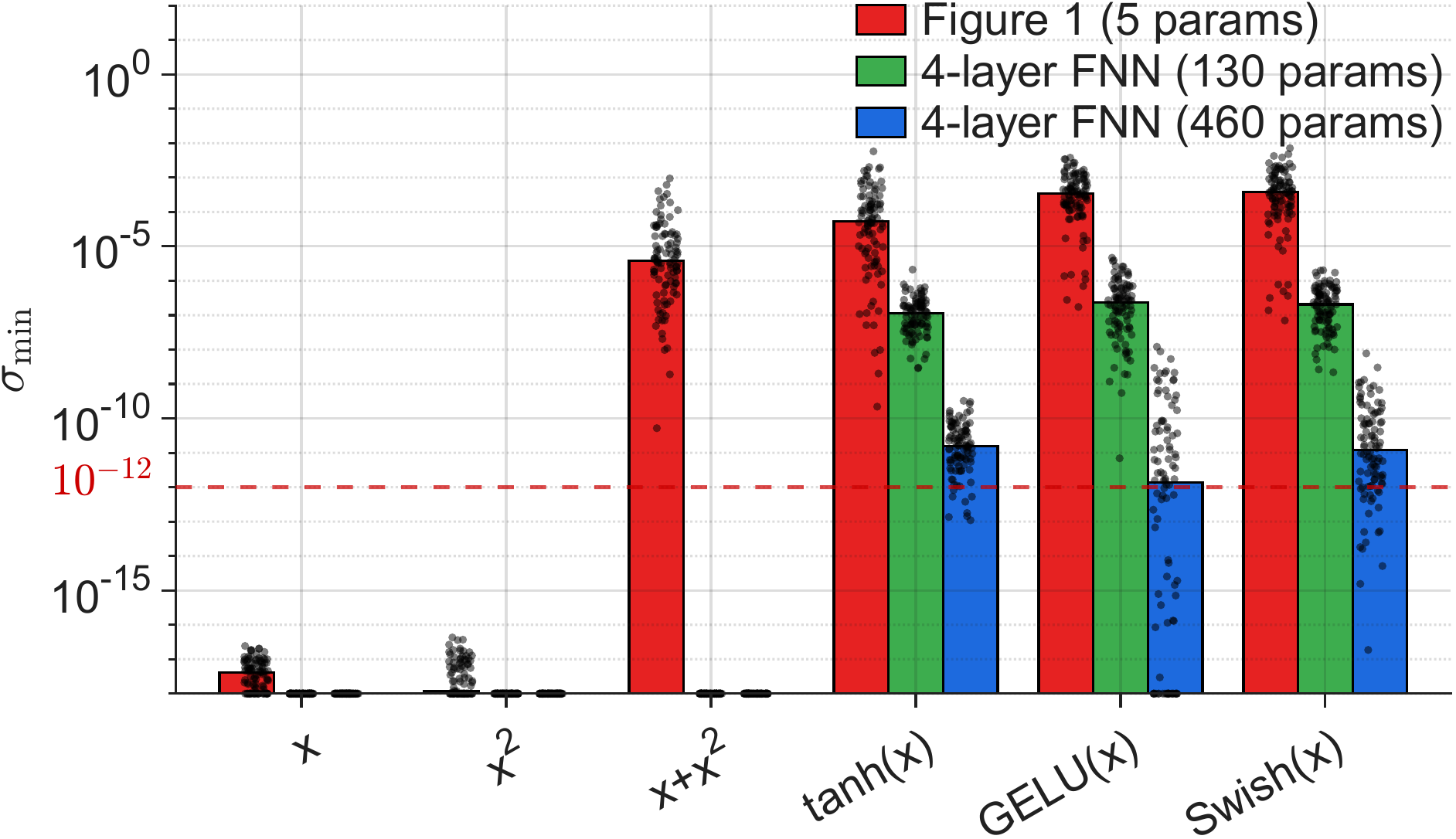}
				%{sim/sim1_2}
				\caption{\black Minimum singular values (median) for $100$  realizations of $W$ for the dynamics in Fig. \ref{fig:ex1} and two FNNs (see Ex. \ref{ex:1})}%\caption{Minimum singular value across different activation functions for $50$ random realizations of $W$ for the dynamics in Figure \ref{fig:ex1}, and two (homogeneous) $4$-layer perceptron with $2-10-10-1$ and $2-20-20-1$ neurons in each layer, respectively.}
				\label{fig:sim}
			\end{figure}

		\addtocounter{example}{-2}
		{\black \begin{example}[cont'd]
		Consider the dynamics in Figure \ref{fig:ex1} with $f_i=f$ for all $i$ and let Assumption $1$ hold. Its condensation graph has a source component $C_1=\{1,2\}$ linked to a sink component $C_2=\{3,4\}$ (see Figure \ref{fig:cond_graph}). Hence, by Theorem \ref{th:1}, the network is locally identifiable for almost all functions in $\mc F_Z$ and almost all weights consistent with the topology when node $2$ is excited and node $3$ is measured. 
		
		To numerically validate this result, we study the local injectivity of the measurement $F_3^{{\black t}}$ in \eqref{eq:ex1} w.r.t. $W$ for various activation functions $f \in \mathcal{F}_Z$ and $100$ random realizations of $W$. For a fixed weight matrix $W$, we evaluate the Jacobian of $F_3^{t}$ with respect to $W$ over $1000$ random input pairs $(u_2^{{\black t}-4},u_2^{{\black t}-2})$ drawn from $[-10, 10]^2$. We then assess the rank of the normalized Jacobian by computing its minimum singular value (a nonzero value yields full-rankness, which guarantees local injectivity). As shown in red in Figure \ref{fig:sim}, the nonlinearities $f(x) \in \{x+x^2, \tanh(x), \mathrm{GELU}(x), \mathrm{Swish}(x)\}$ consistently yield strictly positive minimum singular values, whereas $f(x) = x$ and $f(x) = x^2$ fail to produce full-rank Jacobians, supporting our prior observation that they lack local identifiability. We extend this experiment to fully connected $4$-layer FNNs with widths $2$-$10$-$10$-$1$ and $2$-$20$-$20$-$1$ (shown in green and black in Figure \ref{fig:sim}, respectively). For these wider networks, the numerical behavior changes for $f(x) = x+x^2$, which now fails to produce a full-rank Jacobian. This aligns with the theoretical expectation that the polynomial degree required for identifiability increases with network width. For the remaining nonlinearities, the minimum singular values decrease as width increases (approaching $10^{-12}$), reflecting progressive numerical indistinguishability. Overall, these results empirically support the generic local identifiability established in Theorem \ref{th:1}, while highlighting the gap between theoretical identifiability and numerical identification. %, further discussed in Section \ref{sec:ass}.

			\end{example}}

	The main ideas behind the proof are depicted in the following example.
	\begin{example}[cont'd]
		{\black Consider the dynamics in Figure \ref{fig:path} and let $f\in \mc F_Z$. As observed previously, $F_1^{{\black t}}$ is a composition of analytic functions and is therefore itself analytic. Hence, there exists $r>0$ such that, for all $x \in B_r(0)$,
		\be\label{eq:expF}{\black F_1^t}(x) = A_1(W, a_1)x + A_2(W, a_1, a_2)x^2 +
		\sum_{k>2}A_k(W, a)x^k,\ee
		where $A_k(W,a)$ denotes the $k$-th Maclaurin coefficient of the composite function. %\footnote{To distinguish the intermediate coefficients from the final ones, we denote the Maclaurin coefficients of the composite function by $A_k$.}
	 It depends on $W$ and $a$; e.g., according to \eqref{eq:8},
			$A_1(W,a) = a_1^2 w_{21} {\black w_{12}}$
			and
			$A_2(W,a) = a_2a_1w_{21}^2{\black w_{12}}(1+a_1{\black w_{12}}).$}
		{\black By \eqref{eq:expF}, $F_1^t\equiv\tilde F_1^t$
		if and only if $A_k(W,a)=A_k(\tilde W,a)$} for all $k$. In the following, we prove that, when $f(0)=0$, every coefficient $A_k$ of the measurement $\black F_1^t$ is a polynomial function of the weights $W$ and the first $k$ Maclaurin coefficients $(a_1,\dots,a_k)$ of $f$.
		{\black This allows us to restrict the analysis to a finite number of coefficients indexed by a set $I$ and study local identifiability in terms of the local injectivity in $W$ of the \emph{finite-dimensional} mapping
			\begin{equation}
				A|_I : \mathbb{R}^{|\mathcal{E}|} \times \mathbb{R}^{M} \to \mathbb{R}^{|I|},
			\end{equation}
			which maps the weights $W$ and a finite number $M$ of Maclaurin coefficients of $f$ to the coefficients of $F_1^{{\black t}}$ indexed by $I$.}
		Local injectivity of $A|_I$ with respect to $W$ can be studied by computing the rank of its Jacobian with respect to $W$, denoted by $\nabla_w A|_I(W,a)$. For instance, for $I=\{1,2\}$, %we obtain
	%	\begin{equation*}
	%		\det(\nabla_w A|_I(W,a)) = -a_1^3a_2w_{21}^2{\black w_{12}}.
	%	\end{equation*}
%		By imposing 
%{\black we find the same result as before by imposing } 
{\black we find $\det(\nabla_w A|_I(W,a))=-a_1^3a_2w_{21}^2{\black w_{12}}\ne0$ %, we obtain : the network is generically locally identifiable 
	for every $f$ whose first two Maclaurin coefficients satisfy $a_1\ne0$ and $a_2\ne0$ (as in \eqref{eq:syst_eq}).}
		More generally, we prove that either $\rank(\nabla_w A|_I(W,a))=|\mathcal E|$ for almost all $W$ and almost all $a$, or this equality holds for no $W$ and no $a$. {\black Consequently, to establish generic local identifiability for a given digraph, it is sufficient to exhibit a single activation function and a single weight matrix for which the corresponding network is locally identifiable. In our proof, we use the activation function $f(x)=e^x-1$ as such a witness.}
	\end{example}

			\section{Proof of Theorem \ref{th:1}}\label{sec:proof}
			
			 We begin with a result (proved in Appendix) that allows us to reduce the analysis to %single-{\black excited node} 
			 %and single-{\black measured node} setting 
			 {\black one excited node and one measured node} when the nonlinearity vanishes in zero. {\black Throughout, we always let Assumption \ref{ass1} hold when $\mc G$ contains cycles. }
			
			\begin{lemma}\label{lemma:1s1s}
				Let $\mc G$ be a {\black digraph, with  condensation graph $\mc C(\mc G)$. {\black Assume that one node in each source component of $\mc C(\mc G)$  is excited and one node in each sink component of $\mc C(\mc G)$ is measured}} %. Let Assumption \ref{ass1} hold} 
				and let $f$ be a given nonlinearity satisfying $f(0)=0$.  If, for every excited node $e \in \mc N^e$ and sink $m \in  \mc N^m$, the subgraph	$\mc G^{e\rightarrow m}$ induced by all nodes on a path from $e$ to $m$ is (locally) identifiable, then $\mc G$ is (locally) identifiable.
	\end{lemma}
	
	{\black By Lemma \ref{lemma:1s1s}, proving identifiability for digraphs whose condensation graphs have a single source and a single sink is sufficient to extend the result to all digraphs. Notably, this condition is sufficient but not} necessary, since subgraphs associated to different
	sources and sinks may share common edges that could be identifiable in one subgraph but not in another.  However, since our result holds for {\black all digraphs having condensation graph} with one source and one sink,		this assumption is not restrictive. {\black Therefore, we assume w.l.o.g. that the condensation graph of $\mathcal{G}$ has exactly one source and one sink. Similarly, since the input-output map is required to agree for all input histories, we consider one input coordinate at a time, setting all other inputs to zero. Different coordinates may identify different subsets of weights; it suffices that every weight is identifiable from at least one coordinate. We then simplify the notation for the sink's measurement at time ${\black t}$, $F^{{\black t}}_m$, to $F$, which is now a function from $\R$ to $\R$.} %{\black Nevertheless, because our framework successfully establishes identifiability for all single-source, single-sink digraphs, this structural assumption loses no generality.}
	
%	According to Lemma \ref{lemma:1s1s}, we can assume, w.l.o.g., that the	{\black digraph} $\mc G$ {\black is such that its condensation graph} has one source and one sink. If we prove the result	for this class of graphs, it will naturally extend to all digraphs.		We remark that identifiability with {\black the condensation graph having} with one source and one sink is	sufficient for identifiability with multiple sources and multiple	sinks but not necessary, since subgraphs associated to different sources and sinks may share common edges that could be	identifiable in one subgraph but not in another. However, since	our result holds for {\black all digraphs having condensation graph} with one source and one sink,	this assumption is not restrictive.

			The proof of Theorem \ref{th:1} is now divided in three main steps.
			
			1) In Section \ref{ss:gen_W}, we prove that, for $f$ in $\mathcal{F}_{Z}$, local identifiability is a generic property in {the original parameter space} $W$, that is, it either holds for almost all $W$ or for no $W$. This implies that if local identifiability holds for one %{\black original configuration} 
			$W^{*}$ consistent with the graph {\black with edge set $\mc E$} then local identifiability holds generically for all $W$.
			
			2) In Section \ref{ss:gen_f}, we further prove that $W$-generic local identifiability is a generic property for $f$ in $\mathcal{F}_{Z}$. By combining the two results, we obtain that local identifiability either holds for almost all $W$ for almost all functions $f$ or for no $W$ and no $f$. It is then sufficient to show local identifiability for one $f^{*}$ in $\mathcal{F}_{Z}$ and one $W^{*}$ consistent with $\mathcal{E}$ to obtain $f$-$W$-generic local identifiability of the network $\mathcal{G}$ in the class $\mathcal{F}_{Z}$.
			
			3) Section \ref{ss:exp} is devoted to {\black exhibiting a witness pair $(W^*, f^*)$ for every possible digraph. Indeed, by points 1) and 2), establishing that a digraph $\mc G$ is locally identifiable for at least one such pair yields local identifiability for almost all weights and almost all functions in $\mc F_Z$. In our proof, we choose the witness function as $f^{*}(x)=e^{x}-1$. }%{\black exhibit one witness for every possible digraph. In our proof, this witness is $f^{*}(x)=e^{x}-1$. Then, for every digraph $\mc G$, there exists at least one witness $f^*$ for which the network is identifiable, }thus leading to the main result.  }
			%prove local identifiability {\black in general digraphs for}  $f^{*}(x)=e^{x}-1$. {\black Then, for every digraph $\mc G$, there exists at least one witness $f^*$ for which the network is identifiable, }thus leading to the main result. 

			\subsection{Genericity in the weight matrix}\label{ss:gen_W}
			We will prove that, for a fixed $f$ in $\mathcal{F}_{Z}$, local identifiability is a $W$-generic property. Recall that for $f$ analytic with Maclaurin series coefficients $a=\{a_{k}\}_{k\in\mathbb{N}}$, the measurement $F$ is analytic and %for $x$ sufficiently small 
			%it must take the form
			{\black $\exists r>0$ such that for all $x\in \mc B_r(0)$: }
			\begin{equation}
				F(x) = \sum_{k} A_{k}(W,a)x^{k} \label{eq:12}
			\end{equation}
			for some {\black Maclaurin} coefficients $\black A_{k}(W,a)$ that depend on the weights $W$ and the coefficients $a$. 
			%for some {\black Maclaurin} coefficients
			%\be\label{eq:A_B_map}
			%\black A_{k}(W,a):={B_{k}(\varphi(W),a)}
			%\ee that depend on the weights $W$ {\black (via the linear mapping $\varphi$ in \eqref{eq:phi})} and the coefficients $a$. 
			Let us denote with $A(\cdot,\cdot)$ the operator that maps the weights $W$ and the Maclaurin coefficient of $f$ to the Maclaurin coefficients of $F$. Then, a network $\mathcal{G}$ is locally identifiable for $f$ in $C^{\omega}$ if {\black $\exists\epsilon>0$ such that} for all $\tilde{W}$ {\black consistent with $\mc E$ satisfying $\|\tilde W-W\|<\epsilon$} %consistent with the graph %sufficiently close to $W$, it holds
			\begin{equation}
			A(W,a) = A(\tilde{W},a) \Rightarrow W = \tilde{W}\,. \label{eq:13}
			\end{equation}
			Notice that $A$ is an operator as it maps $W$ and the series of coefficients of $f$ to the series of coefficients of $F$. So {\black the left side of} \eqref{eq:13} is in fact an infinite number of equations, given by the coefficients of $F$ and $\tilde{F}$, and a finite number of variables, that is, the weights in $W$ and $\tilde{W}$. We start by formulating the local identifiability question in terms of local injectivity.
			\begin{definition}
			A function $g:\mathbb{R}^{N}\rightarrow\mathbb{R}^{M}$ is locally injective at $x$ if there exists $\epsilon>0$ such that, for all $\tilde{x}$ satisfying $||x-\tilde{x}||<\epsilon$, 
			%\begin{equation*}
			{\black 	$g(x) = g(\tilde{x})$ implies $x = \tilde{x}\,.$}
			%\end{equation*}
			\end{definition}
			
			To leverage local injectivity in our analysis, we need to consider finite restrictions of the operator $A$ in \eqref{eq:13}. For any finite subset $I\subset\mathbb{N}$ with $|I|=M>0$, we denote with $A|_{I}$ the restriction of $A$ on $I$. We now aim to establish the connection between local identifiability and local injectivity of $A|_{I}$. We first need the following intermediate result, which is proved in Appendix. This result will also be used in the proof of Proposition \ref{pr:genf}.
			
			\begin{lemma}\label{lemma:pol}
			Consider a {\black digraph having condensation graph} with one source and one sink and let $f$ in $\mathcal{F}_{Z}$ with Maclaurin series coefficients $a=\{a_{k}\}_{k\in\mathbb{N}_{+}}$. Then, each coefficient in \eqref{eq:12} is a polynomial function in the weights $W$ and the first $k$ coefficients $(a_{1},\dots,a_{k})$ and does not depend on the higher-order coefficients $\{a_{k^{\prime}}\}_{k^{\prime}>k}$, that is,
			\begin{equation}
				A_{k}(W,a) = A_{k}(W,a_{1},\dots,a_{k}), \label{eq:14}
			\end{equation}
			where $A_{k}:\mathbb{R}^{|\mathcal{E}|}\times\mathbb{R}^{k}\rightarrow\mathbb{R}$ is a polynomial function.
			\end{lemma}
			
			\begin{remark}\label{rem:f0}
			The conclusion of Lemma \ref{lemma:pol} is no longer true when $f(0)\ne0$. %Consider a path graph with 3 nodes 
			{\black Consider the dynamics in Example \ref{ex:path} where the measurement of node $1$ is given by \eqref{eq:F1}. L}et $f$ in $\mathcal{F}$ have Maclaurin series coefficients $a=\{a_{k}\}_{k\in\mathbb{N}}${, \black and suppose that there exists $\rho>0$ such that $w_{32}f(w_{21}x)\in\mathcal B_r(0)$ for all $x\in\mathcal B_\rho(0)$, where $r$ is the radius of convergence of the Maclaurin series of $f$. Then, for such $x$,}
			\begin{align*}
				&F^{{\black t}}_{\black 1}(x) = f(w_{{\black 12}}f(w_{21}x)) = \sum_{k^{\prime}}a_{k^{\prime}}\big(w_{{\black 12}}\sum_{k}a_{k}w_{21}^{k}x^{k}\big)^{k^{\prime}}= \\
				& \sum_{k}a_{k}(a_{0}w_{{\black 12}})^{k} + \big(a_{1}w_{{\black 12}}w_{21}\sum_{k>0}ka_{k}(a_{0}w_{{\black 12}})^{k-1}\big)x + o(x)\,.
			\end{align*}
			Observe that the first coefficient of $\black F_1^t$ is given by $A_{0}(W,a)=\sum_{k}a_{k}(a_{0}w_{32})^{k}$ and depends on all the series of coefficients of $f$, and so does $A_{1}(W,a)$. Furthermore, both $A_{0}$ and $A_{1}$ are analytic in $W$ and $a$.
			\end{remark}
			
			Thanks to Lemma \ref{lemma:pol}, we can establish a link between local identifiability of the network and local injectivity of $A|_{I}$. This is made possible because $A(W,a)$ is polynomial in $W$, allowing us to apply results from algebraic geometry.
			
			\begin{proposition}\label{pr:hilb}
			Consider a {\black digraph} $\mathcal{G}$ {\black having condensation graph} with one source and one sink and let $f$ in $\mathcal{F}_{Z}$. Then, %{\black under Assumption \ref{ass1}},
			 the following three facts are equivalent:
			\begin{enumerate}
			\item[(i)] the network is (locally) identifiable in $W$ by exciting {\black one node in the source component of $\mc C(\mc G)$ and measuring one node in the sink component of $\mc C(\mc G)$;}%the source and measuring the sink;
\item[(ii)]   {\black $\exists\epsilon>0$ such that \eqref{eq:13} holds true for any $\tilde{W}$ consistent with the graph satisfying $\|W-\tilde{W}\|<\epsilon$;}
\item[(iii)]  there exists a finite subset $I\subset\mathbb{N}$ such that $A|_{I}(W,a)$ is (locally) injective in $W$.
			\end{enumerate}
			\end{proposition}
			
			\begin{proof}
			For $f$ analytic, the equivalence between $(i)$ and $(ii)$ is immediate by definition. Moreover, condition $(iii)$ directly implies $(ii)$ as $A|_{I}(W,a)$ is a subset of the conditions in \eqref{eq:13}. Hence we just need to prove {\black $(ii)\Rightarrow(iii)$}. 
			%$(ii)\Rightarrow (iii)$.
			By Lemma \ref{lemma:pol}, each coefficient $A_{k}$ is polynomial in $W$. Therefore, \eqref{eq:13} is a system of an infinite number of polynomial equations. The implication $(ii)\Rightarrow(iii)$ is then true for Hilbert's basis theorem \cite{hilbert1890ueber}, from which it follows that a locus-set of a collection of polynomials is the locus of finitely many polynomials. Then, the solution set of the infinite system in \eqref{eq:13} is determined by finitely many equation, yielding to $(iii)$.
			\end{proof}
			
			According to Proposition \ref{pr:hilb}, the network is (locally) identifiable if and only if $\exists I$, $I$ finite, such that $A|_{I}$ is (locally) injective in $W$. We then focus on $A|_{I}$, for $I$ finite, and we study under which conditions it is local injective in $W$. In the following Lemma, we prove that, for every given $f$ in $C^{\omega}$ and $I$ finite, local injectivity of $A|_{I}$ is a generic property in $W$, i.e., it either holds for almost all weights or for no weights. Also, we show a direct link between local injectivity of $A|_{I}$ in $W$ and full-rankness of its Jacobian, denoted with $\nabla_{w}A|_{I}(W,a)$.
			
			\begin{lemma}\label{lemma:gen}
			Let $\mathcal{G}$ be a {\black digraph having condensation graph} with one source and one sink and let $f$ in $C^{\omega}$ with Maclaurin series coefficients $a=\{a_{k}\}_{k\black \in \N}$. For every $I\subset\mathbb{N}$, $I$ finite, exactly one of the two following holds:
			\begin{itemize}
			\item[(i)] $A|_{I}$ is $W$-generically locally injective in $W$ and $\text{rank } \nabla_{w}A|_{I}(W,a)=|\mathcal{E}|$ for almost all $W$;
			\item[(ii)] there is no $W$ for which $A|_{I}$ is locally injective and $\text{rank } \nabla_{w}A|_{I}(W,a)=|\mathcal{E}|$.
			\end{itemize}
			\end{lemma}
			
			The proof of Lemma \ref{lemma:gen} can be found in Appendix and is based on the results in \cite{Legat2020, Legat2024}. If we combine Proposition \ref{pr:hilb} and Lemma \ref{lemma:gen}, we obtain that local identifiability is a generic property in $W$, that is, it either holds for almost all weights $W$ or for no $W$. This leads to the following proposition.
			
			\begin{proposition}\label{pr:gen_W}
			Consider a {\black digraph having condensation graph $\mc C(\mc G)$} with one source and one sink and let $f$ in $\mathcal{F}_{Z}$. Exactly one of the following holds:
			\begin{itemize}
			\item[(i)] the network is $W$-generically locally identifiable by exciting {\black one node in} the source {\black component} and measuring {\black one node in} the sink {\black component of $\mc C(\mc G)$};
			\item[(ii)] the network is locally identifiable by exciting {\black one node in} the source {\black component} and measuring {\black one node in} the sink {\black component of $\mc C(\mc G)$} for \emph{no $W$}.
			\end{itemize}
			\end{proposition}
			\begin{proof}
			From Lemma \ref{lemma:gen}, it follows that exactly one of the following holds: either $A|_{I}$ is locally injective for no $W$ for all $I$ finite, or there exists $I$ finite such that $A|_{I}$ is locally injective for almost all $W$. If we combine this with Proposition \ref{pr:hilb}, we obtain the result.
			\end{proof}
			
			\begin{remark}
			According to Proposition \ref{pr:gen_W}, for every given $f$ in {\black $\mc F_Z$}, if the network is locally identifiable for one weight matrix $W^{*}$ consistent with the graph, then, the network is locally identifiable for almost all $W$. 			%{\black Because the image of $\varphi$ has measure zero in the codomain $\mathbb{R}^{|\mc E_H|}$, evaluating genericity directly on unfolded matrices $W_H$ would fail. However, because we strictly parameterize by the original edges, the finite variables are precisely the weights in $W$ and  $\tilde W$.}
			\end{remark}
			
			\subsection{Genericity in the functions}\label{ss:gen_f}
			We now prove that local identifiability by exciting sources and measuring sinks is a $f$-generic property in the class $\mathcal{F}_{Z}$, that is, it holds for almost all functions (in the sense of Definition \ref{def:genf}) or for no functions.
			
			\begin{proposition}\label{pr:genf}
			Consider a {\black digraph having condensation graph $\mc C(\mc G)$}  with one source and one sink, %that are excited and measured, respectively. 
			{\black and assume that one node in the source component of $\mc C(\mc G)$  is excited and one node in the sink component of $\mc C(\mc G)$ is measured.}
			Then, exactly one of the following holds:
			\begin{itemize}
			\item[(i)] {\black $\mc G$} is $f$-$W$-generically {\black locally} identifiable in $\mathcal{F}_{Z}$;
		\item[(ii)]  {\black $\mc G$} is {\black locally} identifiable for no $W$ and no $f$ in $\mathcal{F}_{Z}$.
		\end{itemize}
			\end{proposition}
			\begin{proof}
			We remark that $(i)$ and $(ii)$ are mutually exclusive. Then, we can prove the statement with the following reasoning: we assume that $(ii)$ does not hold and prove that then $(i)$ holds, i.e., we assume that the network is locally identifiable for some $W^{*}$ and $f^{*}$ and we show that this implies that the network is locally identifiable for almost all $W$ and almost all $f$. Let $a^{*}$ denote the Maclaurin coefficients of $f^{*}$. By Proposition \ref{pr:hilb}, if the network is locally identifiable for some $W^{*}$, then there exists a finite subset $I\subset\mathbb{N}$ such that $A|_{I}$ is locally injective in $W^{*}$. By Lemma \ref{lemma:gen}, this implies that $A|_{I}$ is $W$-generically locally injective in $W$ and $\text{rank } \nabla_{w}A|_{I}(W,a^{*})=|\mathcal{E}|$ for almost all $W$. In particular, we have that $\text{rank } \nabla_{w}A|_{I}(W^{*},a^{*})=|\mathcal{E}|$. We now want to study the genericity in both $W$ and $a$. Without loss of generality, let us assume that $I=\{1,\dots,M\}$, for some $M>0$. By Lemma \ref{lemma:pol}, we have that $A|_{I}$ depends only on the first $M$ coefficients of $f$ and so we can write $\text{rank } \nabla_{w}A|_{I}(W,a) = \text{rank } \nabla_{w}A|_{I}(W,v)$ for some vector $v$ in $\mathbb{R}^{M}$ that denotes the first $M$ coefficients of $a$. Then, by the genericity of the full-rankness of the Jacobian of an analytic function (see Lemma \ref{lemma:pol}.2 in \cite{Legat2024}), we have that either $(i)$ $\text{rank } \nabla_{w}A|_{I}(W,v)=|\mathcal{E}|$ for almost all $W$ consistent with $\mathcal{E}$ and $v$ in $\mathbb{R}^{M}$ or (ii) $\text{rank } \nabla_{w}A|_{I}(W,v)=|\mathcal{E}|$ for no $W$ and no $v$. Recall that we have $\text{rank } \nabla_{w}A|_{I}(W^{*},v^{*})=|\mathcal{E}|$, for $v^{*}=a^{*}|_{I}$. Then, it must be case $(i)$, that is, $\text{rank } \nabla_{w}A|_{I}(W,v)=|\mathcal{E}|$ for almost all $W$ for almost all $v$ in $\mathbb{R}^{M}$. By Lemma \ref{lemma:gen}, we further have $A|_{I}$ is $W$-generically locally injective in $W$ for almost all $v\in\mathbb{R}^{M}$. According to Proposition \ref{pr:hilb} and Definition \ref{def:genf}, this implies that the network is $W$-generically locally identifiable for almost all functions in $\mathcal{F}_{Z}$. This concludes the proof.
			\end{proof}
			
			According to Proposition \ref{pr:genf}, if the network is locally identifiable for one weight matrix $W^{*}$ consistent with the graph for one $f^{*}$ in $\mathcal{F}_{Z}$ then, the network is locally identifiable for almost all $W$ for almost all $f$.
			
			\subsection{Exponential {\black witness}}\label{ss:exp}	
			Following Proposition \ref{pr:genf}, we now {\black exhibit} one {\black witness activation} function for which any {\black digraph} $\mathcal{G}$ %{\black whose condensation graph contains} one source and one sink {\black component} 
			is $W$-generically locally identifiable {\black by exciting one node in the source component and one node in the sink component of its condensation graph}. This {\black allows us to conclude} that any such %such DAG
			{\black digraph } $\mathcal{G}$ is locally identifiable for almost all $f$ for almost all $W$. %			Specifically, c
			%{\mc C}onsider a {\black digraph} $\mathcal{G}$ {where a node}   $e=1$ {in the source component of $\mc C(\mc G)$ is excited} and one {\black node} $m$ in the sink {\black component} is measured. 
		{\black Let $e=1$ be the excited node %in the source component 
			and $m$ be the measured node. Let} % in the sink component. Let}
			\begin{equation}
			{\black f^*}(x) = e^{x}-1, \label{eq:15}
			\end{equation} {\black and consider a time instant ${\black t}>0$ large enough for all the weights to have an effect on the measurement $F_m^t$ and let all excitation be zero except for one, denoted with $x$.  }
			{\black Following the unfolding procedure in Section \ref{ss:io_map} (see \eqref{eq:6})} and by \eqref{eq:15}, we obtain ${\black F^{t}_{1}}(x)=x$ and, for all $i$ in $\mathcal{N}\backslash\{1\}$,
			\begin{equation*}
			 {\black F^{t}_{i} }(x)   
			{\black = \exp\big(\sum_{j}{\black (\varphi(W))_{ij}}F_{j}(x)\big)-1\,,}
			\end{equation*}
		{\black with $\varphi$ as in \eqref{eq:phi}.} Our identifiability analysis is based on the asymptotic behavior of the functions $F^{{\black t}}_{i}(x)$ as $x\rightarrow+\infty$. Due to the exponential nature of $\black f^*$, we can distinguish nodes and paths based on their asymptotic influence, allowing us to define a hierarchy of stronger paths and nodes.
		
		{\black The unfolding procedure
			replicates each original weight  $w_{ij}$ across all corresponding edges $(i_k, j_{k+1})$ in the unfolded graph. Hence, $\varphi$ is injective and the original weight matrix $W$ can be uniquely recovered from its unfolded representation $\varphi(W)$. Therefore, it suffices to establish identifiability of the unfolded weight matrix $W_H$ on $\operatorname{Im}(\varphi)$. In fact, we prove the stronger statement
			\begin{equation}
				F^t(x,W_H)=F^t(x,\tilde W_H)
				\quad\Longrightarrow\quad
				W_H=\tilde W_H,
				\label{eq:unfold_id}
			\end{equation}
			for arbitrary unfolded weight matrices $W_H$ and $\tilde W_H$ in the class considered below. Equation \eqref{eq:unfold_id} immediately implies identifiability when $W_H,\tilde W_H\in\operatorname{Im}(\varphi)$ and, by injectivity of $\varphi$, identifiability of the original weight matrix $W$.}
			
			As we shall see, our reasoning holds when the weights are all positive and they are different for edges outgoing from the same node. {\black Since we have established that identifiability is a generic property, restricting the parameter space to this subset is sufficient.\footnote{\black In fact, we could restrict it even further to a single witness, since, by genericity, exhibiting one identifiable instance implies identifiability for almost all parameter choices.}}
			Therefore, in the following, we will consider weight matrices $W$ belonging to the set
			\begin{align*}
			\mathcal{W}_{+}^\text{\black o} := \{ W \in \mathbb{R}_{+}^{n\times n}
			 \mid   &\forall i\in\mathcal{N}, w_{j_{1}i} \ne  w_{j_{2}i} \\
			&\forall j_{1}, j_{2}\in\mathcal{N}_{i}^{+}, j_{1} \ne j_{2} \}
			\end{align*} {\black
			In the unfolded graph $H_{m}^{{\black t}}(\mathcal{G})$, two distinct edges outgoing from a node copy $i_{k}$ to copies $j_{k-1}$ and $l_{k-1}$ correspond to two distinct original edges $(i, j)$ and $(i, l)$ with weights $w_{ji}$ and $w_{li}$, respectively. Requiring $w_{j_1i} \ne w_{j_2i}$ in the original space $\mathcal{W}_{+}^\text{o}$ ensures that $w^H_{j_1i} \ne w^H_{j_2i}$ for all corresponding edges in the unfolded graph $H_{m}^{t}(\mathcal{G})$, preserving the  property  for the exponential witness analysis.}
			
			The remainder of this section is organized as follows: we first define the relation of stronger paths and stronger nodes, then we derive the asymptotic properties of $F_{i}$ as $x$ grows to infinity and finally we exploit these properties to show the local identifiability of {\black the unfolded graph $\mathcal{H}_m^t(\mc G)$} for $W_{\black H}$ in $\black \mc W_+:=\varphi(\mathcal{W}_{+}^\text{o})$ and %$f(x)=e^{x}-1$. 
			{\black $f^*$ as in \eqref{eq:15}. In the rest of the proof, to simplify notation, we will refer to the unfolded DAG as $\mc H$ and to its weights as $W_H$ in $\mc W_+$, keeping the weight-sharing constraints of $\text{Im}(\varphi)$ implicit. }
			
			\subsubsection{Stronger paths and stronger nodes}
			Let us begin by defining an order relation on paths. Stronger paths will dominate others when the nonlinearity at the node level is the exponential function. Recall that $\mathcal{P}^{1\rightarrow k}$ denotes the set of paths from 1 to $k$. We denote with
			{\black $\mathcal{P}({\black \mc H}) := \bigcup_{k}\mathcal{P}^{1\rightarrow k}$} the set of all paths outgoing from 1.
			
			\begin{definition}
			Consider a %graph ${\black \mc H}=(\mathcal{N},\mathcal{E},W)$
			{\black DAG $\mc H=(\mc N_H, \mc E_H, W_H )$} with a single source and $W_{\black H}\in\mathcal{W}_{+}$. Let $P_{1}, P_{2}$ in $\mathcal{P}({\black \mc H})$. Then, we say that $P_{1}$ is stronger than $P_{2}$, and we denote $P_{1}\succ^{p}P_{2}$, if $P_{1}\ne P_{2}$ and either
			$(i)$ $|P_{1}|>|P_{2}|$ or
			$(ii)$ $|P_{1}|=|P_{2}|$ and $w_{P_{1}(h)P_{1}(h-1)} > w_{P_{2}(h)P_{2}(h-1)}$
			for $h=\min\{k \mid P_{1}(k)\ne P_{2}(k)\}$.
			\end{definition}
			
			\begin{lemma}\label{lemma:sp}
			The stronger-path relation  $\succ^{p}$ defines a strict total order\footnote{We recall that a binary relation $>$ on a set $X$ is called a \emph{strict total order} if it satisfies the following properties for all $a,b,c \in X$: (i) \emph{Irreflexivity}: $\neg(a>a)$, (ii) \emph{Transitivity}: If $a>b$ and $b>c$, then $a>c$, and (iii) {\black \emph{Totality}}: Exactly one of the following holds: $a>b$, $a=b$, or $b>a$.} on $\mathcal{P}({\black \mc H})$.
			\end{lemma}
			\begin{proof}		$(i)$ \textit{Irreflexivity.} Follows immediately from the definition of $\succ^{p}$.
				$(ii)$ \textit{Transitivity}. Let $P_{1}, P_{2}$ and $P_{3}$ in $\mathcal{P}({\black \mc H})$ be such that $P_{1} \succ^{p} P_{2}$ and $P_{2} \succ^{p} P_{3}$. Then, if $|P_{1}| > |P_{2}|$ or $|P_{2}| > |P_{3}|$ we have $|P_{1}| > |P_{3}|$. Otherwise, if $|P_{1}| = |P_{2}| = |P_{3}|$, {\black let $h = \min\{k \mid P_{1}(k) \ne P_{2}(k) \vee P_{2}(k) \ne P_{3}(k)\}$. Because all three paths are identical prior to step $h$, the weight inequalities at step $h$ must satisfy one of the following mutually exclusive cases: 
					%Then, either
					\begin{equation*}
						\ba
					&	w_{P_{1}(h)P_{1}(h-1)}^{\black H}> w_{P_{2}(h)P_{2}(h-1)}^{\black H} > w_{P_{3}(h)P_{3}(h-1)}^{\black H}\,, \quad \text{or}\\
					&w_{P_{1}(h)P_{1}(h-1)}^{\black H}> w_{P_{2}(h)P_{2}(h-1)}^{\black H} = w_{P_{3}(h)P_{3}(h-1)}^{\black H}\,, \quad \text{or}\\
					&	w_{P_{1}(h)P_{1}(h-1)}^{\black H}= w_{P_{2}(h)P_{2}(h-1)}^{\black H} > w_{P_{3}(h)P_{3}(h-1)}^{\black H}.
						\ea
					\end{equation*}
					In all cases, $w_{P_{1}(h)P_{1}(h-1)}^{\black H} > w_{P_{3}(h)P_{3}(h-1)}^{\black H}$.  Since $P_1(k)=P_3(k)$ for all $k < h$, it follows that $h = \min\{k \mid P_{1}(k) \ne P_{3}(k)\}$, implying $P_{1} \succ^{p} P_{3}$.
				%	In all cases, $P_{1} \succ^{p} P_{3}$.
				$(iii)$  \textit{Totality}. Let $P_{1}$ and $P_{2}$ in $\mathcal{P}({\black \mc H})$ be such that $P_{1} \ne P_{2}$. We assume w.l.o.g. that $|P_{1}| \ge |P_{2}|$. Then, either $|P_{1}| > |P_{2}|$, and thus $P_{1} \succ^{p} P_{2}$, or $|P_{1}| = |P_{2}| = \ell$. In this latter case, $\exists h \in \{1, \dots, \ell\}$ such that $h = \min\{k \mid P_{1}(k) \ne P_{2}(k)\}$. Therefore, either
			%	\begin{equation*}
			%		\ba
			%		w_{P_{1}(h)P_{1}(h-1)}^{\black H} > w_{P_{2}(h)P_{2}(h-1)}^{\black H} \quad&\Rightarrow\quad P_{1} \succ^{p} P_{2}\,,\quad \text{or}\\
			%		w_{P_{2}(h)P_{2}(h-1)}^{\black H} > w_{P_{1}(h)P_{1}(h-1)}^{\black H}\quad&\Rightarrow\quad P_{2} \succ^{p} P_{1}\,.
			%		\ea
			%	\end{equation*}
			$w_{P_{1}(h)P_{1}(h-1)}^{\black H} > w_{P_{2}(h)P_{2}(h-1)}^{\black H}$ which implies $ P_{1} \succ^{p} P_{2}$, or
				$w_{P_{2}(h)P_{2}(h-1)}^{\black H} > w_{P_{1}(h)P_{1}(h-1)}^{\black H}$ which implies $P_{2} \succ^{p} P_{1}$.}				\end{proof}
			
			Since $\succ^{p}$ is a strict total order and $\mathcal{P}({\black \mc H})$ is a finite set, the strongest path $\black 			P_{*}({\black \mc H}) := \max_{\succ^{p}} \mathcal{P}({\black \mc H})$
			 of a graph ${\black \mc H}$
			%\begin{equation*}
			%\end{equation*}
			is well-defined. We now define the notion of stronger nodes.
			
			\begin{definition}\label{def:sp}
			We say that a node $i$ in $\mathcal{N}_{\black H}$ is \emph{stronger} that $j$, and we denote $i\succ^{n}j$, if
$	\black			P_{*}({\black \mc H}^{1\rightarrow i}) \succ^{p} P_{*}({\black \mc H}^{1\rightarrow j})$.
			\end{definition}
			
			We have that also $\succ^{n}$ defines a strict total order, as proved in the following lemma.
			
			\begin{lemma}\label{lemma:sto}
			For every {\black DAG} ${\black \mc H=(\mc N_H, \mc E_H, W_H)}$ with $W_{\black H}\in\mathcal{W}_{+}$ with a single source, the relation $\succ^{n}$ is a strict total order on the set $\mathcal{N}_{\black H}$.
			\end{lemma}
			\begin{proof}
			Follows from the fact that $\succ^{p}$ is a strict total order and $P_{*}({\black \mc H}^{1\rightarrow i})=P_{*}({\black \mc H}^{1\rightarrow j})$ if and only if $i=j$.
			\end{proof}
			
			Stronger paths and stronger nodes satisfy the following properties.
			
			\begin{lemma}\label{lemma:prop}
			For a node $i$ in $\mathcal{N}_{\black H}$, let $P_{i}=P_{*}({\black \mc H}^{1\rightarrow i})$, $l_{i}=|P_{i}|$ and $i^{\prime}=P_{i}(l_{i}-1)$. Then,
			$(i)$ $i^{\prime}\succ^{n}k$ for all $k$ in $\mathcal{N}_{i}^{-}$, $k\ne i^{\prime}$.
			Let $j$ be such that $i\succ^{n}j$. Then,
			$(ii)$ ${\black w^H_{ii^{\prime}}}>{\black w^{ H}_{ji^{\prime}}}$, and
			$(iii)$ for $j^{\prime}=P_{j}(l_{j}-1)$, either $i^{\prime}=j^{\prime}$ or $i^{\prime}\succ^{n}j^{\prime}$.
			\end{lemma}
			\begin{proof}
			$(i)$ Assume by absurd that $k\succ^{n}i^{\prime}$. Then, $P_{i}^{\prime}=(P_{k},i)$ is such that $P_{i}^{\prime}\in\mathcal{P}^{1\rightarrow i}$ and $P_{i}^{\prime}\succ^{p}P_{i}$, thus contradicting the assumption that $P_{i}$ is the strongest path to $i$.
			$(ii)$ Observe that, by definition of $\mathcal{W}_{+}$, we have ${\black w^H_{ii^{\prime}}}\ne {\black w^H_{ji^{\prime}}}$ since $i\ne j$. Assume now by absurd that ${\black w^H_{ii^{\prime}}<w^H_{ji^{\prime}}}$. Then, $P_{j}=(P_{i^{\prime}},j)$ satisfies $P_{j}\in\mathcal{P}^{1\rightarrow j}$ and $P_{j}\succ^{p}P_{i}$. This leads to $j\succ^{n}i$, which is a contradiction.
			$(iii)$ If $|P_{i}|>|P_{j}|$, we have $|P_{i^{\prime}}|=|P_{i}|-1>|P_{j^{\prime}}|=|P_{j}|-1$ and therefore $i^{\prime}\succ^{n}j^{\prime}$. Otherwise, if $|P_{i}|=|P_{j}|$, we have that $h=\min\{k \mid P_{i}(k)\ne P_{j}(k)\}$ either satisfies $h=l_{i}$ and therefore $i^{\prime}=j^{\prime}$, or $h<l_{i}$ and therefore $i^{\prime}\succ^{n}j^{\prime}$.
			\end{proof}
			
			\subsubsection{Asymptotic behavior}
			We now present the first key asymptotic result justifying the order introduced on nodes.
			
			\begin{proposition}\label{prop:asymp1}
			Consider a graph ${\black \mc H}$ with a single source and $W_{\black H}\in\mathcal{W}_{+}$, and let $f(x)=e^{x}-1$. Then, for $i, j$ in $\mathcal{N}_{\black H}$,
			\begin{equation}
				i\succ^{n}j \Rightarrow \lim_{x\rightarrow+\infty}\frac{F_{j}(x)}{F_{i}(x)}=0. \label{eq:16}
			\end{equation}
			\end{proposition}
			\begin{proof}
			For a node $k$ in $\mathcal{N}$, let $P_{k}=P_{*}({\black \mc H}^{1\rightarrow k})$ denote the strongest path from the source 1 to $k$, and let $l_{k}=|P_{k}|$ denote its length. We proceed by induction on $l_{i}$.
			
			\underline{Base case}: $l_{i}=1$. Then $P_{i}=(1,i)$, so node $i$ is directly connected to the source. There are two possibilities for $j$:\\
			\emph{Case 1:} $l_{j}=0$. Then, $j=1$, i.e., $F_{j}(x)=x$ and
			\begin{equation*}
				\lim_{x\rightarrow+\infty}\frac{F_{j}(x)}{F_{i}(x)} = \lim_{x\rightarrow+\infty}\frac{x}{\exp({\black w^{H}_{i1}}x)-1} = 0.
			\end{equation*}
			\emph{Case 2}: $l_{j}=1$. Then, $P_{j}=(1,j)$ and
			$$			%\begin{align*}
				\lim_{x\rightarrow+\infty}\frac{F_{j}(x)}{F_{i}(x)} =%&= 
				\lim_{x\rightarrow+\infty}\frac{\exp(w_{j1}^{\black H}x)-1}{\exp(w_{i1}^{\black H}x)-1}\,.
			$$	
			 %\\&= \lim_{x\rightarrow+\infty}\exp((w_{j1}^{\black H}-w_{i1}^{\black H})x)\,.
			%\end{align*}
			Since $i\succ^{n}j$, Lemma \ref{lemma:prop}$(ii)$ implies $w_{i1}^{\black H}>w_{j1}^{\black H}$. Hence, $\lim_{x\rightarrow+\infty}\frac{F_{j}(x)}{F_{i}(x)}=0$.
			
			\underline{Inductive step}: We now consider node $i$ with $l_{i}>1$ and we let $j$ be such that $i\succ^{n}j$. We assume that the result holds for all nodes $k^{\prime}\succ^{n}k$, with $l_{k^{\prime}}\le l_{i}-1$. Let $i^{\prime}=P_{i}(l_{i}-1)$ be the predecessor of $i$ along the strongest path. Then $l_{i^{\prime}}=l_{i}-1$, so by the induction hypothesis, for any $k$ such that $i^{\prime}\succ^{n}k$,
			\begin{equation}
				\lim_{x\rightarrow+\infty}\frac{F_{k}(x)}{F_{i^{\prime}}(x)}=0. \label{eq:17}
			\end{equation}
			From Lemma \ref{lemma:prop}$(i)$, we have that $i^{\prime}\succ^{n}k$ for all $k\in\mathcal{N}_{i}^{-}$. Analogously, letting $j^{\prime}=P_{j}(l_{j}-1)$, we find $j^{\prime}\succ^{n}k$ for all $k\in\mathcal{N}_{j}^{-}$. From Lemma \ref{lemma:prop}$(iii)$, we further have $i^{\prime}\succ^{n}j^{\prime}$. Combining all the conditions together, we find $i^{\prime}\succ^{n}k$ for all $k\in\mathcal{N}_{i}^{-}\cup\mathcal{N}_{j}^{-}$, and the inductive assumption \eqref{eq:17} holds for all such $k$.
			Now consider the limit:
			\begin{align*}
				\lim_{x\rightarrow+\infty}\frac{F_{j}(x)}{F_{i}(x)} &= \lim_{x\rightarrow+\infty}\frac{\exp(\sum_{k}w_{jk}^{\black H}F_{k}(x))-1}{\exp(\sum_{k}w_{ik}^{\black H}F_{k}(x))-1} \\
				&= \lim_{x\rightarrow+\infty}\exp\big(\sum_{k}(w_{jk}^{\black H}-w_{ik}^{\black H})F_{k}(x)\big)\,.
			\end{align*}
			We rewrite this sum by factoring out $F_{i^{\prime}}(x)$, thus obtaining
			\begin{align*}
					\lim_{x\rightarrow+\infty}\frac{F_{j}(x)}{F_{i}(x)} &=\lim_{x\rightarrow+\infty} \exp\Big(((w_{ji^{\prime}}-w_{ii^{\prime}})F_{i^{\prime}}(x)  \\
				&\qquad \big(1+\sum_{k\in\mathcal{N}_{i}^{-}\cup\mathcal{N}_{j}^{-}}\frac{(w_{jk}^{\black H}-w_{ik}^{\black H})F_{k}(x)}{(w_{ji^{\prime}}^{\black H}-w_{ii^{\prime}}^{\black H})F_{i^{\prime}}(x)}\big)\Big)\,.
			\end{align*}
			The induction hypothesis implies that each term in the sum goes to zero as $x\rightarrow+\infty$, so the expression in brackets converges to 1. Moreover, since $i\succ^{n}j$, Lemma \ref{lemma:prop}$(ii)$ ensures $w_{ji^{\prime}}<w_{ii^{\prime}}$, so the overall exponent tends to $-\infty$, and the exponential tends to zero. Hence, $\lim_{x\rightarrow+\infty}\frac{F_{j}(x)}{F_{i}(x)}=0$. %This completes the proof.
			\end{proof}
			
			We now state the second key asymptotic result, which is based on Proposition \ref{prop:asymp1} and serves as the main technical foundation for the identifiability analysis in the following section.
			
			\begin{proposition}\label{prop:asymp2}
			Consider a {\black DAG $\mc H$} with a single source and $W_{\black H}\in\mathcal{W}_{+}$, and let $f(x)=e^{x}-1$. Then, for $i, j$ in $\mathcal{N}$,
			\begin{equation}
				\lim_{x\rightarrow+\infty}\frac{\log(F_{i}(x))-\sum_{k\succ^{n}j}w_{ik}^{\black H}F_{k}(x)}{F_{j}(x)} = w_{ij}^{\black H}. \label{eq:18}
			\end{equation}
			\end{proposition}
			\begin{proof}
			Let us denote with $k\prec^{n}j$ when $j\succ^{n}k$. Then,
			\begin{align*}
				&\lim_{x\rightarrow+\infty}\frac{\log(F_{i}(x))-\sum_{k\succ^{n}j}w_{ik}^{\black H}F_{k}(x)}{F_{j}(x)} \\
				&= \lim_{x\rightarrow+\infty}\frac{\sum_{k}w_{ik}^{\black H}F_{k}(x)-\sum_{k\succ^{n}j}w_{ik}^{\black H}F_{k}(x)}{F_{j}(x)} \\
				&= \lim_{x\rightarrow+\infty}\frac{w_{ij}^{\black H}F_{j}(x)+\sum_{k\prec^{n}j}w_{ik}^{\black H}F_{k}(x)}{F_{j}(x)} \\
				&= \lim_{x\rightarrow+\infty}\frac{w_{ij}^{\black H}F_{j}(x)\big(1+\sum_{k\prec^{n}j}\frac{w_{ik}^{\black H}F_{k}(x)}{w_{ij}^{\black H}F_{j}(x)}\big)}{F_{j}(x)}\,.
			\end{align*}
			By Proposition \ref{prop:asymp1}, each term in the sum goes to zero as $x\rightarrow+\infty$ and the expression in brackets converges to 1. Hence, \eqref{eq:18} holds.
			\end{proof}
			
			\subsubsection{Identifiability}
			%Consider a DAG ${\black \mc H}$ with a single source and single sink. In this section, we aim to 
			{\black We now} conclude our reasoning by showing that {\black every DAG ${\black \mc H}$ with a single source and a single sink} is locally identifiable for $f(x)=e^{x}-1$ by exciting the source and measuring the sink, that is, %$\tilde{W}_{\black H}$ sufficiently close to $W_{\black H}$ then
			{\black there exists an $\epsilon>0$ such that} for all $\tilde{W}_H$ {\black consistent with the graph satisfying $\|\tilde W_H-W_H\|<\epsilon$ it holds that} $F \equiv \tilde{F}$ {\black implies} $W_{\black H} = \tilde{W}_{\black H}$.
			%\begin{equation*}
			%F \equiv \tilde{F} \Rightarrow W_{\black H} = \tilde{W}_{\black H}.
			%\end{equation*}
			In order to apply the previous asymptotic results, we need the order among the nodes to be preserved.
			
			\begin{definition}
			Let ${\black \mc H}=(\mathcal{N}_{\black H},\mathcal{E}_{\black H},W_{\black H})$ be a DAG with a single source, and let $W_{\black H} \in \mathcal{W}_{+}$. We say that $\tilde{W}_{\black H}$ in $\mathcal{W}_{+}$ preserves the order of $W_{\black H}$ if $i\succ^{n}j$ in ${\black \mc H}$ if and only if $i\succ^{n}j$ in $\tilde{{\black \mc H}}=(\mathcal{N},\mathcal{E},\tilde{W})$.
			\end{definition}
			
			We now employ Proposition \ref{prop:asymp2} and, in particular, equation \eqref{eq:18} to derive Lemma 8, which provides the core reasoning behind the proof of local identifiability. Recall that, for a node $i$ in $\mathcal{N}$, $W^{1\rightarrow i}$ denotes the weight matrix of the subgraph ${\black \mc H}^{1\rightarrow i}$ induced by the set of nodes that belong to paths that go from the source 1 to the node $i$.
			
			\begin{lemma}
			Consider a DAG ${\black \mc H}=(\mathcal{N}_{\black H},\mathcal{E}_{\black H},W_{\black H})$ with a single source and a single sink, and a weight matrix $W$ in $\mathcal{W}_{+}$. Let $f(x)=e^{x}-1$ and let $\tilde{W}_{\black H}\in\mathcal{W}_{+}$ be consistent with the graph and preserve the order of $W_{\black H}$. Take any node $i\ne1$. If the output of the node $i$ satisfy
			\begin{equation}
				\lim_{x\rightarrow+\infty}\frac{F_{i}(x)}{\tilde{F}_{i}(x)} = c_{i} \ne 0 \label{eq:19}
			\end{equation}
			while the induced weight matrices on the subgraph from the source to $i$ differ, that is,
			%\begin{equation}
			$\black 	W^{1\rightarrow i}_{\black H} \ne \tilde{W}^{1\rightarrow i}_{\black H}, \label{eq:20}$
			%\end{equation}
			then, there exists an in-neighbor $j$ in $\mathcal{N}_{i}^{-}$ such that the same holds, that is,
			\begin{equation}
				W^{1\rightarrow j}_{\black H} \ne \tilde{W}^{1\rightarrow j}_{\black H}, \quad \lim_{x\rightarrow+\infty}\frac{F_{j}(x)}{\tilde{F}_{j}(x)} = c_{j} \ne 0. \label{eq:21}
			\end{equation}
			\end{lemma}
			
			In words, if a node $i$ has the same asymptotic output under $W$ and $\tilde{W}$ but the weights on the paths to $i$ differ, this behavior must already appear one step earlier, at some in-neighbor $j$. This lemma {\black will allow} the backward recursion used in Proposition \ref{prop:exp_rec}: if $W_{\black H} \ne \tilde{W}_{\black H}$, Lemma 8 lets us propagate the two properties until we eventually reach the out-neighbors of the source of the graph, whose asymptotic output under different $W$ and $\tilde{W}$ cannot match, contradicting $F \equiv \tilde{F}$.
			
			\begin{proof}
			We start by proving
			\begin{equation}
				\mc T = \{k \in \mathcal{N}_{i}^{-} \mid W^{1\rightarrow k}_{\black H} \ne \tilde{W}^{1\rightarrow k}_{\black H}\} \ne \emptyset\,.  \label{eq:22}
			\end{equation}
			Assume by contradiction that $W^{1\rightarrow k}_{\black H} = \tilde{W}^{1\rightarrow k}_{\black H}$ for all $k$ in $\mathcal{N}_{i}^{-}$. Since by assumption $W^{1\rightarrow i}_{\black H} \ne \tilde{W}^{1\rightarrow i}_{\black H}$, this implies that there must exist an edge $(i, k)$ such that $w_{ik}^{\black H} \ne \tilde{w}_{ik}^{\black H}$. Let us then define
			\begin{equation}
				j^{*} = \max_{\succ^{n}}\{k \in \mathcal{N}_{i}^{-} \mid w_{ik}^{\black H} \ne \tilde{w}_{ik}^{\black H}\}. \label{eq:23}
			\end{equation}
			According to \eqref{eq:18} in Proposition \ref{prop:asymp2}, it holds that
			\begin{align}
				w_{ij^{*}}^{\black H} &= \lim_{x\rightarrow+\infty}\frac{\log(F_{i}(x)) - \sum_{k\succ^{n}j^{*}} w_{ik}^{\black H}F_{k}(x)}{F_{j^{*}}(x)} \notag \\
				&= \lim_{x\rightarrow+\infty} \frac{\log\big(\frac{F_{i}(x)}{\tilde{F}_{i}(x)}\tilde{F}_{i}(x)\big) - \sum_{k\succ^{n}j^{*}} w_{ik}^{\black H}F_{k}(x)}{F_{j^{*}}(x)} \notag \\
				&\stackrel{(1)}{=} \lim_{x\rightarrow+\infty} \frac{\log(c_{i}) + \log(\tilde{F}_{i}(x)) - \sum_{k\succ^{n}j^{*}} w_{ik}^{\black H}F_{k}(x)}{F_{j^{*}}(x)} \notag \\
				&= \lim_{x\rightarrow+\infty} \frac{\log(\tilde{F}_{i}(x)) - \sum_{k\succ^{n}j^{*}} w_{ik}^{\black H}F_{k}(x)}{F_{j^{*}}(x)}, \label{eq:24}
			\end{align}
			where (1) holds by the assumption on the asymptotic behavior of $F_{i}$ and $\tilde{F}_{i}$ in \eqref{eq:19}. By definition of $j^{*}$, we have that $w_{ik}^{\black H} = \tilde{w}_{ik}^{\black H}$ for all $k \succ^{n} j^{*}$. Also, we assumed that $W^{1\rightarrow k}_{\black H} = \tilde{W}^{1\rightarrow k}_{\black H}$ for all $k$ in $\mathcal{N}_{i}^{-}$ and this, by definition of $F_{k}$, implies that $F_{k} = \tilde{F}_{k}$ for all $k$ in $\mathcal{N}_{i}^{-}$. Using these equalities, we find
			\begin{align}
				{\black w_{ij^{*}}^{H}} &= \lim_{x\rightarrow+\infty}\frac{\log(\tilde{F}_{i}(x)) - \sum_{k\succ^{n}j^{*}} w_{ik}^{\black H}F_{k}(x)}{F_{j^{*}}(x)}  \label{eq:25}
				 \\
				&= \lim_{x\rightarrow+\infty}\frac{\log(\tilde{F}_{i}(x)) - \sum_{k\succ^{n}j^{*}} \tilde{w}_{ik}^{\black H}\tilde{F}_{k}(x)}{\tilde{F}_{j^{*}}(x)} = \tilde{w}_{ij^{*}}^{\black H},\notag			\end{align}
			which contradicts the definition of $j^{*}$ in \eqref{eq:23}, i.e., $w_{ij^{*}}^{\black H} \ne \tilde{w}_{ij^{*}}^{\black H}$. Hence, the set $\mc T$ in \eqref{eq:22} must be nonempty. Let $j = \max_{\succ^{n}}\mc T$. By definition of $\mc T$, the left-hand condition in \eqref{eq:21} is automatically satisfied for this $j$. We now show that the right-hand condition in \eqref{eq:21} also holds for this same $j$. By definition of $j$ and the standing assumption on $F_{i}$,
			\begin{align}
				&\lim_{x\rightarrow+\infty} \big(\log(F_{i}(x)) - \sum_{k\succ^{n}j}w_{ik}^{\black H}F_{k}(x)\big) \notag \\
				&= \lim_{x\rightarrow+\infty} \big(\log(c_{i}) + \log(\tilde{F}_{i}(x)) - \sum_{k\succ^{n}j}w_{ik}^{\black H}\tilde{F}_{k}(x)\big). \label{eq:26}
			\end{align}
			To proceed, we show that $w_{ik}^{\black H} = \tilde{w}_{ik}^{\black H}$ for all $k \succ^{n} j$. If this already holds for all $k$ in $\mathcal{N}_{i}^{-}$, there is nothing to prove. Otherwise, let $j^{*}$ be defined as in \eqref{eq:23} and assume by absurd that $j^{*} \succ^{n} j$. By applying the same reasoning as before, we find
			\begin{align*}
				w_{ij^{*}} ^{\black H}&= \lim_{x\rightarrow+\infty}\frac{\log(F_{i}(x)) - \sum_{k\succ^{n}j^{*}} w_{ik}^{\black H}F_{k}(x)}{F_{j^{*}}(x)} \\
				&= \lim_{x\rightarrow+\infty}\frac{\log(\tilde{F}_{i}(x)) - \sum_{k\succ^{n}j^{*}} \tilde{w}_{ik}^{\black H}\tilde{F}_{k}(x)}{\tilde{F}_{j^{*}}(x)}= \tilde{w}_{ij^{*}}^{\black H},
			\end{align*}
			again contradicting the definition of $j^{*}$ in \eqref{eq:23}. Therefore either $j^{*} = j$ or $j \succ^{n} j^{*}$. In both cases, $w_{ik}^{\black H} = \tilde{w}_{ik}^{\black H}$ for all $k \succ^{n} j$. Then, starting from \eqref{eq:26}, we find
			\begin{align*}
				&1 = \lim_{x\rightarrow+\infty} \frac{\log(F_{i}(x)) - \sum_{k\succ^{n}j} w_{ik}^{\black H}F_{k}(x)}{\log(F_{i}(x)) - \sum_{k\succ^{n}j} w_{ik}^{\black H}F_{k}(x)} \\
				&= \lim_{x\rightarrow+\infty} \frac{\log(F_{i}(x)) - \sum_{k\succ^{n}j} w_{ik}^{\black H}F_{k}(x)}{\log(\tilde{F}_{i}(x)) - \sum_{k\succ^{n}j} \tilde{w}_{ik}^{\black H}\tilde{F}_{k}(x)} =\lim_{x\rightarrow+\infty} \frac{w_{ij}^{\black H}F_{j}(x)}{\tilde w_{ij}^{\black H}\tilde{F}_{j}(x)},
			\end{align*}
			which implies
			%\begin{equation*}
			$\black 	\lim_{x\rightarrow+\infty} \frac{F_{j}(x)}{\tilde{F}_{j}(x)} = \frac{\tilde{w}_{ij}^{\black H}}{w_{ij}^{\black H}} = c_{j} \ne 0.$
			%\end{equation*}
			This concludes the proof.
			\end{proof}
			
			We can now prove local identifiability of ${\black \mc H}$.
			
			\begin{proposition}\label{prop:exp_rec}
			Let ${\black\mathcal{H}}$ be a {\black DAG} with a single source, a single sink and $W_{\black H} \in \mathcal{W}_{+}$ and let $f(x) = e^{x} - 1$. Then, ${\black \mathcal{H}}$ is locally identifiable by exciting the source and measuring the sink.
			\end{proposition}
			\begin{proof}
			We prove the proposition by contradiction. Let $\tilde{W}_{\black H} \ne W_{\black H} $ in $\mathcal{W}_{+}$ be any weight matrix consistent with ${\black \mathcal{H}}$ and %sufficiently close to $W_{\black H} $
		{\black $\epsilon >0$ be such that  every $\tilde{W}_{\black H} $ satisfying $\|\tilde{W}_{\black H} -W_{\black H} \|<\epsilon$,}  preserves the order of $W_{\black H} $. %Assume that the input–output maps coincide, that is, $F \equiv \tilde{F}$, even though $W_{\black H}  \ne \tilde{W}_{\black H} $.
		{\black Assume that $F \equiv \tilde{F}$. } We now show that this assumption leads to an absurd. Since $F_{n} \equiv \tilde{F}_{n}$ for the sink $n$, we have
			\begin{equation*}
				\lim_{x\rightarrow+\infty}\frac{F_{n}(x)}{\tilde{F}_{n}(x)} = 1 \ne 0, \quad W^{1\rightarrow n}_{\black H}  \ne \tilde{W}^{1\rightarrow n}_{\black H} ,
			\end{equation*}
			so Lemma 8 applies with $i = n$ and $c_{n} = 1$. Lemma 8 then yields an in-neighbor $j(n) \in \mathcal{N}_{n}^{-}$ such that
			\begin{equation*}
				W^{1\rightarrow j(n)} \ne \tilde{W}^{1\rightarrow j(n)}_{\black H} , \quad \lim_{x\rightarrow+\infty} \frac{F_{j(n)}(x)}{\tilde{F}_{j(n)}(x)} = c_{j(n)} \ne 0.
			\end{equation*}
			Applying Lemma \ref{lemma:prop} again to $i = j(n)$, and iterating, we obtain a sequence of nodes in the DAG. Because the graph is acyclic, this process must reach some out-neighbor $j \in \mathcal{N}_{1}^{+}$ of the source, still satisfying
			\begin{equation*}
				W^{1\rightarrow j}_{\black H}  \ne \tilde{W}^{1\rightarrow j}_{\black H} , \quad \lim_{x\rightarrow+\infty}\frac{F_{j}(x)}{\tilde{F}_{j}(x)} = c_{j} \ne 0.
			\end{equation*}
			On the other hand, for any out-neighbor of the source, $F_{j}(x) = e^{w_{1j}^{\black H}x} - 1$ {\black and} $\tilde{F}_{j}(x) = e^{\tilde{w}_{1j}^{\black H}x} - 1$, so
			\begin{equation*}
				\lim_{x\rightarrow+\infty} \frac{F_{j}(x)}{\tilde{F}_{j}(x)} = \lim_{x\rightarrow+\infty} \frac{e^{w_{1j}^{\black H}x} - 1}{e^{\tilde{w}_{1j}^{\black H}x} - 1} = 
				\begin{cases} 
					0 & \text{if } w_{1j}^{\black H} < \tilde{w}_{1j}^{\black H}, \\
					+\infty & \text{if } w_{1j}^{\black H} > \tilde{w}_{1j}^{\black H}, \\
					1 & \text{if } w_{1j}^{\black H} = \tilde{w}_{1j}^{\black H}.
				\end{cases}
			\end{equation*}
			The limit can be a finite nonzero constant $c_{j}$ only when $w_{1j} = \tilde{w}_{1j}$, contradicting $W^{1\rightarrow j}_{\black H}  \ne \tilde{W}^{1\rightarrow j}_{\black H} $. Hence $W_{\black H}  = \tilde{W}_{\black H} $.
			\end{proof}
			
			We now conclude by combining Proposition \ref{pr:genf} and 6. By Proposition \ref{prop:exp_rec}, for any DAG ${\black \mathcal{H}}$ with one source and one sink, there exist $f^{*}$ in $\mathcal{F}_{Z}$ and $W^{*}_{\black H} $ {\black in $\text{Im}(\varphi)$} %consistent with the graph 
			such that {\black the unfolded graph is locally identifiable. As a consequence, there exist $f^*$ and $W^*$ consistent with the graph such that the original network $\mathcal{G}$ is locally identifiable}. %is locally identifiable. 
			{\black B}y Proposition \ref{pr:genf}, every {\black digraph having condensation graph} with one source and one sink is $f$-$W$-generically locally identifiable in $\mathcal{F}_{Z}$. By Lemma \ref{lemma:1s1s}, the results extend to {\black digraphs} with multiple sources and multiple sinks. This concludes the proof.
			
			\section{Heterogeneous Nonlinearities}\label{sec:het}
			In the previous sections, we assumed that all nodes share a common nonlinearity $f$. We now consider a generalized model where each node $i \in \mathcal{N}$ is associated with a potentially different nonlinear activation function $f_{i}$ as in (1). All previous assumptions remain unchanged: the graph $\mathcal{G}$ is a {\black digraph}, the functions $f_{i}$ are known for all nodes, and the weight matrix $W$ is consistent with the known topology. %{\red If space allows, add example}
			The definition of identifiability in Definition \ref{def:id} naturally extends to this setting. Specifically, consider a graph $\mathcal{G} = (\mathcal{N}, \mathcal{E}, W)$ with sets of excited and measured nodes $\mathcal{N}^{e}$ and $\mathcal{N}^{m}$, respectively, and let $\mathcal{F}$ be a class of functions. We say that a collection of functions $\{f_{i}\}_{i\in\mathcal{N}}$ belongs to $\mathcal{F}^{\mathcal{N}}$ if $f_{i} \in \mathcal{F}$ for all $i$ in $\mathcal{N}$. Then, a network $\mathcal{G}$ with heterogeneous nonlinearities is identifiable in $\mathcal{F}$ if, for any given collection $\{f_{i}\}_{i\in\mathcal{N}}$ in $\mathcal{F}^{\mathcal{N}}$, {\black there exists ${\black t}>0$ such that $F^{{\black t}}(\mathcal{N}^{m}) = \tilde{F}^t(\mathcal{N}^{m})$} implies $W = \tilde{W}$. 
						\addtocounter{example}{+1}
			{\black 	%Consider a path graph with $4$ nodes where we excite the source $1$ and measure the sink $4$, yielding $F^k_4(x) = f_4(w_{43} f_3(w_{32} f_2(w_{31} x)))$, and let \be \label{eq:nonl}f_i(x) = a_{1}^{(i)}x + a_{2}^{(i)}x^2\ee for $i \in \{1, 2, 3\}$. By  solving the system of equations formed by matching polynomial coefficients of $F^k_4$, we find that the network is identifiable for all $f$ in $\mc F=\{ f_i \text{ as in }\eqref{eq:nonl},, a_{1}^{(1)}a_{2}^{(1)}  a_{1}^{(2)}a_{2}^{(2)} a_{1}^{(3)} \neq 0\}$, 					meaning that inner layers must retain non-linear terms ($a_{12}, a_{22} \neq 0$), while the outer layer $f_3$ can be strictly linear ($a_{32} = 0$). 	This relaxation is impossible in the homogeneous setting where, for $f_i(x) = f(x) = a_1 x + a_2 x^2$ for all $i \in \{1, 2, 3\}$, the network is identifiable if and only if  $a_1a_2 \neq 0$.
					\begin{example}\label{ex:multif}
						Consider a path graph with $4$ nodes where we excite the source node $1$ and measure the sink node $4$, yielding $F_4^{{\black t}}(x) = f_4(w_{43} f_3(w_{32} f_2(w_{21} x)))$. Let 
%						\be \label{eq:nonl}
						$f_i(x) = a_{1}^{(i)}x + a_{2}^{(i)}x^2$
%						\ee 
						for $i \in \{1, 2, 3, 4\}$. By solving the system of equations formed by matching polynomial coefficients of $F_4^{{\black t}}$ and $\tilde F_4^{{\black t}}$, we find that the network is identifiable for all collections $\{f_i\}_{i=1}^4$ in $\mathcal{F} = \{ f_i  \mid a_{1}^{(2)}a_{2}^{(2)} a_{1}^{(3)}a_{2}^{(3)} a_{1}^{(4)} \neq 0\}$, i.e., the inner layers must retain non-linear terms ($a_{2}^{(2)}, a_{2}^{(3)} \neq 0$), while $f_4$ can be linear ($a_{2}^{(4)} = 0$). Note that this relaxation is impossible in the homogeneous setting where, for $f_i(x) = f(x) = a_1 x + a_2 x^2$ for all $i$, the network is identifiable if and only if $a_1 a_2 \neq 0$.
					\end{example}}
		
		Similarly, we can extend the definition of local identifiability in Definition \ref{def:loc_id} and generic identifiability in Definition \ref{def:genW}.
		It is slightly more delicate to extend the definition of genericity in the class of functions (Definition \ref{def:genf}). %Let us assume that each nonlinearity $f_{i} \in C^{\omega}$ is analytic, with Maclaurin series coefficients $a^{(i)} = \{a_{k}^{(i)}\}_{k\in\mathbb{N}}$. Let
			{\black Let $f_{i} \in C^{\omega}$ have Maclaurin coefficients $a^{(i)} = \{a_{k}^{(i)}\}_{k\in\mathbb{N}}$. We denote with}
			\begin{equation}
			\textbf{a} = \{a_{k}^{(i)}\}_{k\in\mathbb{N}, i\in\mathcal{N}} \label{eq:27}
			\end{equation}
			 the series of all Maclaurin coefficients. For $M > 0$, we let $I = \{1, \dots, M\}$ and we denote with
			%\begin{equation}
			$\textbf{a}|_{I} = \{a_{k}^{(i)}\}_{k\in I, i\in\mathcal{N}} \label{eq:28}$
			%\end{equation}
			the \emph{first $M$ Maclaurin coefficients} of $\{f_{i}\}_{i\in\mathcal{N}}$. Note that $\textbf{a}|_{I} \in \mathbb{R}^{nM}$, where $n = |\mathcal{N}|$ denotes the number of nodes. %We have the following definition.
			\begin{definition}
			Consider a class of functions $\mathcal{F} \subseteq C^{\omega}$ and a set $\mathcal{N}$ with $n = |\mathcal{N}|$. A property holds for almost all collections of functions $\{f_{i}\}_{i\in\mathcal{N}}$ in $\mathcal{F}^{\mathcal{N}}$ if {\black $\exists M > 0$} such that the property holds for all collections of functions $\{f_{i}\}_{i\in\mathcal{N}}$ in $\mathcal{F}^{\mathcal{N}}$ except possibly a subset of collections whose first $M$ Maclaurin coefficients all lie on a zero measure set in $\mathbb{R}^{nM}$.
			\end{definition}
			
			We then have the following result.
			
			\begin{theorem}\label{th:het}
			Let $\mathcal{G}$ be a %{\black DAG} 
			{\black digraph with condensation graph $\mc C(\mc G)$} and let $\mathcal{F}_{Z}$ be as in \eqref{eq:9}. Then, {\black under Assumption \ref{ass1}}, the network $\mathcal{G}$ with heterogeneous nonlinearities is $W$-generically locally identifiable for almost all collections $\{f_{i}\}_{i\in\mathcal{N}}$ of functions in $\mathcal{F}_{Z}^{\mathcal{N}}$ if and only if {\black at least one node in each source component of $\mc C(\mc G)$  is excited and one node in each sink component of $\mc C(\mc G)$ is measured.} %{\black at least one node is excited in } all the sources {components of $\mc C(\mc G)$} and {\black at least one node is measured in } all the sink {components of $\mc C(\mc G)$} 
			\end{theorem}
			
			\begin{remark}
			The result holds for almost all collections of functions in $\mathcal{F}_{Z}^{\mathcal{N}}$ except possibly a subset of those whose first $M$ coefficients all lie in a measure zero set of $\mathbb{R}^{nM}$. Note that the set containing the first $M$ coefficients of homogeneous nonlinearities, that is,
			\begin{equation*}
				E = \{(a_{k}^{(i)})_{k\in I,i\in\mathcal{N}} \mid a_{k}^{(i)} =a_{k}^{(j)}, \forall i, j \in \mathcal{N}, \forall k\}
			\end{equation*}
			is a measure zero set of $\mathbb{R}^{nM}$. In other words, the set of homogeneous nonlinearities $\mathcal{F}_{Z,h}^{\mathcal{N}} := \{\{f\}_{i\in\mathcal{N}}, f \in \mathcal{F}_{Z}\}$ has measure zero in the set $\mathcal{F}_{Z}^{\mathcal{N}}$. Therefore, unlike what one might at first expect, Theorem \ref{th:1} does not follow from Theorem \ref{th:het} and the two results are distinct.  %{\black For instance, in Example \ref{ex:multif}, Theorem \ref{th:het} does not guarantee that the network is identifiable for homogeneous choices $a_1^{(i)} = a_1$ and $a_2^{(i)} = a_2$ for all $i$; in fact, no such $\textbf{a}$ may exist under its scope. This is instead guaranteed by Theorem \ref{th:1}.}
			\end{remark}
			\begin{proof}
			The proof of Theorem \ref{th:het} follows the same steps as the one of Theorem \ref{th:1}. %First, observe that the input-output map in {\black Section \ref{ss:io_map}} extends naturally to the case where each node has its own nonlinearity $f_{i}$. 
			As in the homogeneous case, {\black w.l.o.g.}, we assume that $\mathcal{G}$ has a {\black condensation graph with a} single source {\black component} and a single sink {\black one}. %We obtain $F_{1}(x) = x$ and
			%\begin{equation*}
			%$\black 	F_{i}(x) = f_{i}\big(\sum_{j}{\black (\varphi(W))_{ij}}F_{j}(x)\big)$,
			%\end{equation*}
			%for all $i$ in $\mathcal{N} \setminus \{1\}$. As before, our goal is to determine conditions under which $F \equiv \tilde{F}$ implies $W = \tilde{W}$. 
			Consider a collection of analytic nonlinearities $\{f_{i}\}_{i\in\mathcal{N}}$ and let $\textbf{a}$ as in \eqref{eq:27} denote the series of all Maclaurin coefficients. Then, %for small enough inputs $x$, 
			the measurement $F$ is analytic and %admits the form
			{\black $\exists r>0$ such that for all $x\in \mc B_r(0)$:}
%			\begin{equation}
$ \black F(x) = \sum_{k} A_{k}(W, \textbf{a})x^{k}$,
%			\end{equation}
			where the coefficients $A_{k}$ depend both on the weight matrix $W$ and the series of all Maclaurin coefficients $\textbf{a}$. %As before, we denote with $A(\cdot, \cdot)$ the operator that maps the weights and the coefficients of the nonlinearities to the output series of $F$. 
			Then, the network $\mathcal{G}$ is locally identifiable with node dependent analytic nonlinearities if {\black there exists $\epsilon>0$ such that}, for all $\tilde{W}$ {\black satisfying $\|\tilde W-W\|<\epsilon$,} the condition
			\begin{equation}
				A(W, \textbf{a}) = A(\tilde{W}, \textbf{a}) \Rightarrow W = \tilde{W} \label{eq:30}
			\end{equation}
			holds. We now extend Lemma \ref{lemma:pol} to the case of heterogeneous nonlinearities. The proof is provided in Appendix.
			
			\begin{lemma}\label{lemma:het}
				Consider a {\black digraph} $\mathcal{G}$ with {\black condensation graph with} one source and one sink {\black and assume that one node in the source component of $\mc C(\mc G)$  is excited and one node in the sink component of $\mc C(\mc G)$ is measured.}
				 Let $\{f_{i}\}_{i\in\mathcal{N}}$ in $\mathcal{F}_{Z}^{\mathcal{N}}$ be a collection of functions with Maclaurin series coefficients $\textbf{a}$ as in \eqref{eq:27}. Then, each coefficient {\black of $F$} %in \eqref{eq:29} 
				 is polynomial in the weights $W$ and the first $k$ coefficients $(a_{1}^{(i)}, \dots,a_{k}^{(i)})$ of each function $f_{i}$ and does not depend on high-order coefficients: %, that is,
				\begin{equation}
					A_{k}(W, \textbf{a}) = A_{k}(W, \textbf{a}|_{I}) \label{eq:31}
				\end{equation}
				where $I = \{1, \dots, k\}$ and $A_{k} : \mathbb{R}^{|\mathcal{E}|} \times \mathbb{R}^{nk} \rightarrow \mathbb{R}$ is a polynomial function.
			\end{lemma}
			
			Thanks to Lemma \ref{lemma:het} and \eqref{eq:30}, we can derive the following result, which is the extended version of Proposition \ref{pr:genf}.
			
			\begin{proposition}[Extension of Proposition \ref{pr:genf}]\label{prop:het}
				Let $\mathcal{G}$ be a {\black digraph having condensation graph} with one source and one sink {\black component}, and let each node $i \in \mathcal{N}$ have a nonlinearity $f_{i} \in \mathcal{F}_{Z}$. Then, exactly one of the following holds:
				\begin{itemize}
				\item[(i)] The network is $W$-generically locally identifiable for almost all collections $\{f_{i}\}_{i\in\mathcal{N}}$ in $\mathcal{F}_{Z}^{\mathcal{N}}$;
				\item[(ii)] The network is locally identifiable for no $W$ and no collection $\{f_{i}\}_{i\in\mathcal{N}}$ in $\mathcal{F}_{Z}^{\mathcal{N}}$.
				\end{itemize}
			\end{proposition}
			
			\begin{proof}
				The proof follows the same steps as the proof of Proposition \ref{pr:gen_W} with minor modifications. We assume that the network is locally identifiable for some $W^{*}$ and a collection $\textbf{f}^{*} = \{f_{i}^{*}\}_{i\in\mathcal{N}}$ and we show that this implies that the network is locally identifiable for almost all $W$ and almost all $\{f_{i}\}_{i\in\mathcal{N}}$. Let $\textbf{a}^{*}$ denote the Maclaurin coefficients of $\textbf{f}^{*}$. By applying Proposition \ref{pr:hilb} to \eqref{eq:30} instead of \eqref{eq:13} and by Lemma \ref{lemma:gen}, we find that there exists a finite subset $I \subset \mathbb{N}$ such that $A|_{I}$ is locally injective in $W^{*}$ and that $\text{rank } \nabla_{w}A|_{I}(W^{*}, \textbf{a}^{*}) = |\mathcal{E}|$. 
				
				We now want to study the genericity in both $W$ and $\textbf{a}$. Without loss of generality, let us assume that $I = \{1, \dots, M\}$, for some $M > 0$. By Lemma \ref{lemma:het}, we have that $A|_{I}$ depends only on the first $M$ coefficients of the collection $\{f_{i}\}_{i\in\mathcal{N}}$ and therefore $\text{rank } \nabla_{w}A|_{I}(W, \textbf{a}) = \text{rank } \nabla_{w}A|_{I}(W, \textbf{v})$ for some $v$ in $\mathbb{R}^{Mn}$. Then, by the same reasoning as in the proof of Proposition \ref{pr:genf}, we find that $\text{rank } \nabla_{w}A|_{I}(W, \textbf{v}) = |\mathcal{E}|$ for almost all $W$ for almost all $\textbf{v}$ in $\mathbb{R}^{nM}$. By Lemma \ref{lemma:gen}, Proposition \ref{pr:hilb} and Definition \ref{def:genf}, this implies that the network is $W$-generically locally identifiable for almost all collections of functions $\{f_{i}\}_{i\in\mathcal{N}}$ in $\mathcal{F}_{Z}^{\mathcal{N}}$. This concludes the proof.
			\end{proof}
			
			Then, Theorem \ref{th:het} is obtained by combining Proposition \ref{prop:het} with Proposition \ref{prop:exp_rec}.
			\end{proof}
\section{Discussion on the Assumptions}\label{sec:ass}

The results of Theorems \ref{th:1} and \ref{th:het} hold for (almost all) analytic functions that vanish at zero, under the assumption that no biases are present in the network. {\black{ Analyticity in a neighborhood of zero is essential because our proof relies on Maclaurin coefficients; hence, non-smooth functions such as ReLU fall outside our scope, whereas smooth analytic counterparts such as GELU and Swish are valid.}
{\black Our second key assumption is} that $f(0)=0$ and the network has no biases. These conditions ensure that every Maclaurin coefficient of the measurement $F$ in \eqref{eq:12} is polynomial in $W$ and depends only on finitely many coefficients of $f$ (see Lemmas \ref{lemma:pol} and \ref{lemma:het}). {\black As illustrated in Remark \ref{rem:f0}, this property can fail when $f(0)\neq0$: even for the two-node graph in Figure \ref{fig:path}, measurement coefficients may depend on infinitely many coefficients of $f$ and be analytic, rather than polynomial, in $W$. The same issue arises with biases, since the network may evaluate $f$ away from zero even when $f(0)=0$\footnote{\black E.g., for $F(x)=f(w_2f(w_1x+b_1)+b_2)$, we have $A_0=f(w_2f(b_1)+b_2)$ and $A_1=w_1w_2f'(b_1)f'(w_2f(b_1)+b_2)$.}.}

The polynomial property of the coefficients is crucial {\black for} Proposition \ref{pr:hilb}, which {\black allows us to} extract finitely many equations from the {\black infinite} system \eqref{eq:13}. {\black By Hilbert's basis theorem, finitely many polynomial equations suffice to describe the solution set of an infinite system of polynomial equations. This equivalence fails when the coefficients $A_k(W,a)$ are analytic but not polynomial in $W$.} For example, consider $\sin(kw)=0$ for all $k\in\mathbb{N}$. The full solution set is $w=0$, whereas any finite subset has a strictly larger solution set. This distinction is critical because our identifiability proof is explicitly established for the exponential function, which is analytic. {\black Without a reduction to finitely many equations, the result for the exponential function might not extend to other analytic functions.}
{\black Extending the analysis beyond $f(0)=0$ and zero biases is left for future work. One possible direction is to determine whether local identifiability for one polynomial activation implies local identifiability for almost all analytic activations. This would reduce the problem to finding a suitable polynomial witness for each graph, which remains open. The results of \cite{Finkel2024, Usevich2025} do not provide such witnesses, since monomial activations yield identifiability only up to rescaling. We conjecture that the desired result may nevertheless hold when the polynomial degree is sufficiently large.}

Beyond these constraints, our framework assumes perfect prior knowledge of the node nonlinearities and graph topology. If the activation function is unknown or contains unknown parameters, rescaling issues may arise, making it structurally impossible to separate $W$ from $f$ when some nodes are neither measured nor excited. Additional measurements or excitations may therefore be required for identifiability as in \cite{Vizuete2024a}. {\black The assumption of known topology $\mc E$ raises a further open question. Generic identifiability of the weights for a given topology does not imply identifiability of the architecture. An overparameterized network may represent the same input-output map as a smaller one by setting additional weights to zero. Such configurations form a measure-zero subset of the parameter space and are therefore excluded by generic identifiability results, even though their identifiability is essential for recovering the underlying architecture. Extending the framework to architecture identifiability would thus require addressing these nongeneric cases and equivalences between different topologies.}

Another key assumption is local identifiability. As discussed in Section \ref{ss:id}, we focus on local identifiability because global identifiability is generally impossible even for simple network structures (see Example \ref{ex:ann}). {\black The results could be strengthened by considering identifiability in equivalence classes. This would, however, require different proof techniques, since our local analysis relies on characterizations through local injectivity and full-rank conditions.} {\black Finally, identifiability should be distinguished from practical identification. Our results establish theoretical uniqueness of the network parameters but do not provide an identification algorithm. A direct implementation of the proof would require extracting Maclaurin coefficients from high-order derivatives of the measurements, which is highly sensitive to noise. Thus, our graph-theoretic result should be viewed as providing minimal theoretical conditions. Robust numerical identification will likely require exciting and measuring a larger subset of nodes, which we leave for future work.}}

\section{Conclusions}\label{sec:conc}

{\black We studied the identifiability of networks with linear edge dynamics and nonlinear node dynamics, given the topology and node nonlinearity, and we proved that digraphs are generically locally identifiable for almost all analytic functions that vanish at zero by exciting only one node in each source component of the condensation graph and measuring one node in each sink component. This result holds for both homogeneous and heterogeneous nonlinearities. It contrasts with prior results showing that network identifiability requires exciting and/or measuring every node for both linear and nonlinear dynamic and extends prior work on ANN identifiability to a broader class of functions and network structures. Future work includes extending the analysis to broader classes of nonlinearities and incorporating offsets, as discussed in Section \ref{sec:ass}.}

			\bibliographystyle{ieeetr}
			\bibliography{bib}
			
			\appendix
		%	\subsection{Proofs}
			\subsubsection*{Proof of Lemma \ref{lemma:1s1s}}\label{app:a}
			{\black Let us fix a time ${\black t}$ and let} {\black $F^{{\black t}}$} denote the input-output map of $\mathcal{G}$ {\black for such ${\black t}$}. For any source $e \in \mathcal{N}^{e}$ and sink $m \in \mathcal{N}^{m}$, let $\black F_{e\rightarrow m}^{ {\black t}}$ denote the input-output map restricted to the subgraph $\mathcal{G}^{e\rightarrow m}$. Assume that $F^{\black {\black t}} \equiv \tilde{F}^{\black t}$ for all inputs. Then, in particular, for the sink $m$, we have $F_{m} ^{\black t}\equiv \tilde{F}_{m}^{\black t}$. Now, consider setting all inputs except $e$ to zero, i.e., $x^{(i)} = 0$ for all $i \ne e$. Since $f(0) = 0$, contributions from all other sources vanish, and we obtain $F_{m}^{\black t}(x) = F_{e\rightarrow m}^{\black t}(x^{(e)})$ and $\tilde{F}_{m}^{\black t}(x) = \tilde{F}_{e\rightarrow m}^{\black t}(x^{(e)})$.
			By assumption, each subgraph ${\black \mc H}^{e\rightarrow m}$ is identifiable. Therefore, {\black there exists a $t^{e\rightarrow m}$ such that}
			\begin{equation*}
				F_{e\rightarrow m}^{\black t^{e\rightarrow m}} \equiv \tilde{F}_{e\rightarrow m}^{\black t^{e\rightarrow m}} \Rightarrow W^{e\rightarrow m} = \tilde{W}^{e\rightarrow m}.
			\end{equation*}
			This holds for all source-sink pairs $(e, m)$. {\black Since} every edge in the network {\black belongs to} one subgraph ${\black \mc H}^{e\rightarrow m}$, all weights in ${\black \mc G}$ are {\black identifiable by considering $t^*=\max_{e,m}t^{e\rightarrow m}$, i.e., {\black $F^{t^*}\equiv\tilde F^{t^*}$} implies $W = \tilde{W}$.}
			
			\subsubsection*{Proof of Lemma \ref{lemma:het} and Lemma \ref{lemma:pol}}
{\black W}e prove Lemma \ref{lemma:het}, while Lemma \ref{lemma:pol} follows as a special case. Let $\mathcal{L}(i, j) = \max\{\vert{}P\vert{} \mid P \in \mathcal{P}^{i\rightarrow j}\}$ denote the length of the longest path from node $i$ to node $j$, and let $L = \mathcal{L}(1, m)$ be the length of the longest path between the source node (node 1) and the sink node $m$. We {\black partition} the node set into layers $\mathcal{N} = \mathcal{N}^{(0)} \cup \dots \cup \mathcal{N}^{(L)}$, {\black where $\mathcal{N}^{(l)} = \{i \in \mathcal{N} \mid \mathcal{L}(1, i) = l\}$.} Note that for any node $i \in \mathcal{N}^{(l)}$ in layer $l$, every in-neighbor $j$ of $i$ belongs to an earlier layer $\mathcal{N}^{(l')}$ with $l' < l$. We prove the statement by strong induction on the layer index $l$.
	
\underline{Base case.} For $l = 0$, {\black	$\mathcal{N}^{(0)} = \{1\}$ contains only the source. Since $F_1(x) = x=%$ and it can be rewritten as $
\sum_{k>0} A_k^{(1)}(W, a_1,\dots,a_k)x^k$ with $A_1^{(1)}=1$ and $A_k^{(1)}=0$ $\forall k>0$, the statement holds for $l=0$.}
	
	\underline{Inductive step.} Consider now a layer $l \in \{1, \dots, L\}$. For every $j \in \mathcal{N}^{(0)} \cup \dots \cup \mathcal{N}^{(l-1)}$, {\black denote the collection of the first $k$ coefficients of $\{f_{j'}\}_{j'\in S^{1\rightarrow j}}$ by $\textbf{a}_k^{1\rightarrow j} := \{a_{k'}^{(j')}\}_{k'\le k, , j'\in S^{1\rightarrow j}}$.	By the inductive hypothesis, for all $j \in \mathcal{N}^{(0)} \cup \dots \cup \mathcal{N}^{(l-1)}$, 	$F_j(x) = \sum_{k>0} A_k^{(j)}(W, \textbf{a}_k^{1\rightarrow j}) x^k$, 
	where each coefficient $A_k^{(j)}$ is polynomial function in the weights $W$ and in $\textbf{a}_k^{1\rightarrow j}$. Then, for all $i \in \mathcal{N}^{(l)}$, we have
$$	\begin{aligned}
		&F_i(x) = f_i( \sum_j \varphi(W)_{ij} F_j(x) ) \\&= f_i( \sum_j \varphi(W)_{ij} ( \sum_{k>0} A_k^{(j)}(W, \textbf{a}_k^{1\rightarrow j}) x^k ) ) = f_i( \sum_{k>0} b_k^{(i)} x^k ),\end{aligned}$$ where $b_k^{(i)} := \sum_j \varphi(W)_{ij} A_k^{(j)}(W, \textbf{a}_k^{1\rightarrow j})$ depend only on $W$ and $\mathbf{a}_k^{1\rightarrow j}$ and $b_0^{(i)}=0$. Using the Taylor expansion $f_i(y) = \sum_{k>0} a_k^{(i)} y^k $ and setting $y = \sum_{k>0} b_k^{(i)} x^k$, we find:\begin{align*}F_i(x) &= \sum_{k=1}^M a_k^{(i)} ( \sum_{k'>0} b_{k'}^{(i)} x^{k'})^k + \sum_{k>M} a_k^{(i)} ( \sum_{k'>0} b_{k'}^{(i)} x^{k'})^k \\&\overset{(1)}{=} \sum_{k=1}^M a_k^{(i)} ( \sum_{k'=1}^M b_{k'}^{(i)} x^{k'}+\sum_{k'>M} b_{k'}^{(i)} x^{k'} )^k + o(x^M)\\&\overset{(2)}{=} \sum_{k=1}^M a_k^{(i)} ( \sum_{k'=1}^M b_{k'}^{(i)} x^{k'} )^k + o(x^M)\,,\end{align*}
		where $(1)$ and $(2)$ hold since $k>0$ (that is, $b_0^{(i)}=0$ and $a_0^{(i)}=0$).	
		%Expanding the finite polynomial $\sum_{k=1}^M a_k^{(i)} \left( \sum_{k'=1}^M b_{k'}^{(i)} x^{k'} \right)^k$ and gathering all terms of degree $k \le M$, all higher-order cross-terms of degree $> M$ are absorbed into $o(x^M)$. 
		In words, all powers of degree $k>M$ give rise to monomials of degree $k>M+1$. If we combine this with the definition of $b_k^{(i)}$, % depend on $W$ and $\textbf{a}_k^{1\rightarrow i}$, 
		we find $F_i(x) = \sum_{k=1}^M A_k^{(i)}(W, \textbf{a}_k^{1\rightarrow i}) x^k + o(x^M)$
	where each coefficient $A_k^{(i)}(W, \mathbf{a}_k^{1\rightarrow i})$ %for $1 \le k \le M$ is a
	is polynomial in $W$ and $\mathbf{a}_k^{1\rightarrow i}$, depending solely on the first $k$ coefficients of $\{f_{j'}\}_{j'\in S^{1\rightarrow i}}$. %Since $M \ge 1$ was arbitrary, the coefficient functions $A_k^{(i)}$ are well-defined for all $k \ge 1$ and independent of $M$. 
	Taking $i = m$ concludes the proof.}
			\subsubsection*{Proof of Lemma \ref{lemma:gen}}
			The proof of Lemma \ref{lemma:gen} is adapted from the proof of Theorem 4.1 in \cite{Legat2024} and makes use of Proposition \ref{pr:genf}.1 in \cite{Legat2024} (also Proposition \ref{pr:genf}.1 in \cite{Legat2020}) and Lemma \ref{lemma:pol}.2 in \cite{Legat2024}. In order to use these results, we need to introduce the notion of coordinate-injectivity {\black \cite{Legat2020, Legat2024}}
%			\begin{definition}[\cite{Legat2020}]
				A function $g : \mathbb{R}^{N} \rightarrow \mathbb{R}^{M}$ is \emph{locally coordinate-injective} for coordinate $e$ at $x \in \mathbb{R}^{N}$ if there exists $\epsilon > 0$ such that, for all $\tilde{x}$ satisfying $\|x - \tilde{x}\| < \epsilon$, there holds
				%\begin{equation*}
				$ \black 	g(x) = g(\tilde{x}) \Rightarrow x_{e} = \tilde{x}_{e}$. %
				%\end{equation*}
	%		\end{definition}
		%	\begin{proof}[Proof of Lemma \ref{lemma:gen}]
				For $a$ and $I$ fixed, let us define $g(W) := A|_{I}(W, a)$. Following \cite{Legat2024}, we denote with $e_{e}$ the standard basis vector filled with zeros except 1 at the $e$-th entry. Then, according to Proposition \ref{pr:genf}.1 in \cite{Legat2024}, since $g : \mathbb{R}^{n\times n} \rightarrow \mathbb{R}^{M}$ is an analytic function where $\mathbb{R}^{n\times n}$ and $\mathbb{R}^{M}$ are smooth manifolds with finite dimension, and $\mathbb{R}^{n\times n}$ has dimension $|\mathcal{E}|$ by the rank-nullity theorem. This implies that either (i) $\text{rank } \nabla g(W) = |\mathcal{E}|$ for almost all $W$, or $(ii)$ $\text{rank } \nabla g(W) = |\mathcal{E}|$ for almost no $W$. On the other hand, from Lemma \ref{lemma:pol}.2 in \cite{Legat2024}, we know that the full rankness of $\nabla g(W)$ is a generic property, that is, it either holds for all $W$ or for no $W$. Therefore, we obtain that either (i) $\text{rank } \nabla g(W) = |\mathcal{E}|$ for almost all $W$, or $(ii)$ $\text{rank } \nabla g(W) < |\mathcal{E}|$ for all $W$. In case $(ii)$, this implies that for each $W$ there is a direction $e$ for which $\ker \nabla g(W) \not\perp e_{e}$, since the rank is non-full. We apply this argument for each $W$, and obtain that in case $(ii)$ $g$ is locally injective for no $W$. This concludes the proof.
\begin{IEEEbiography}
	[{\includegraphics[width=1in,height=1.25in,clip,keepaspectratio]{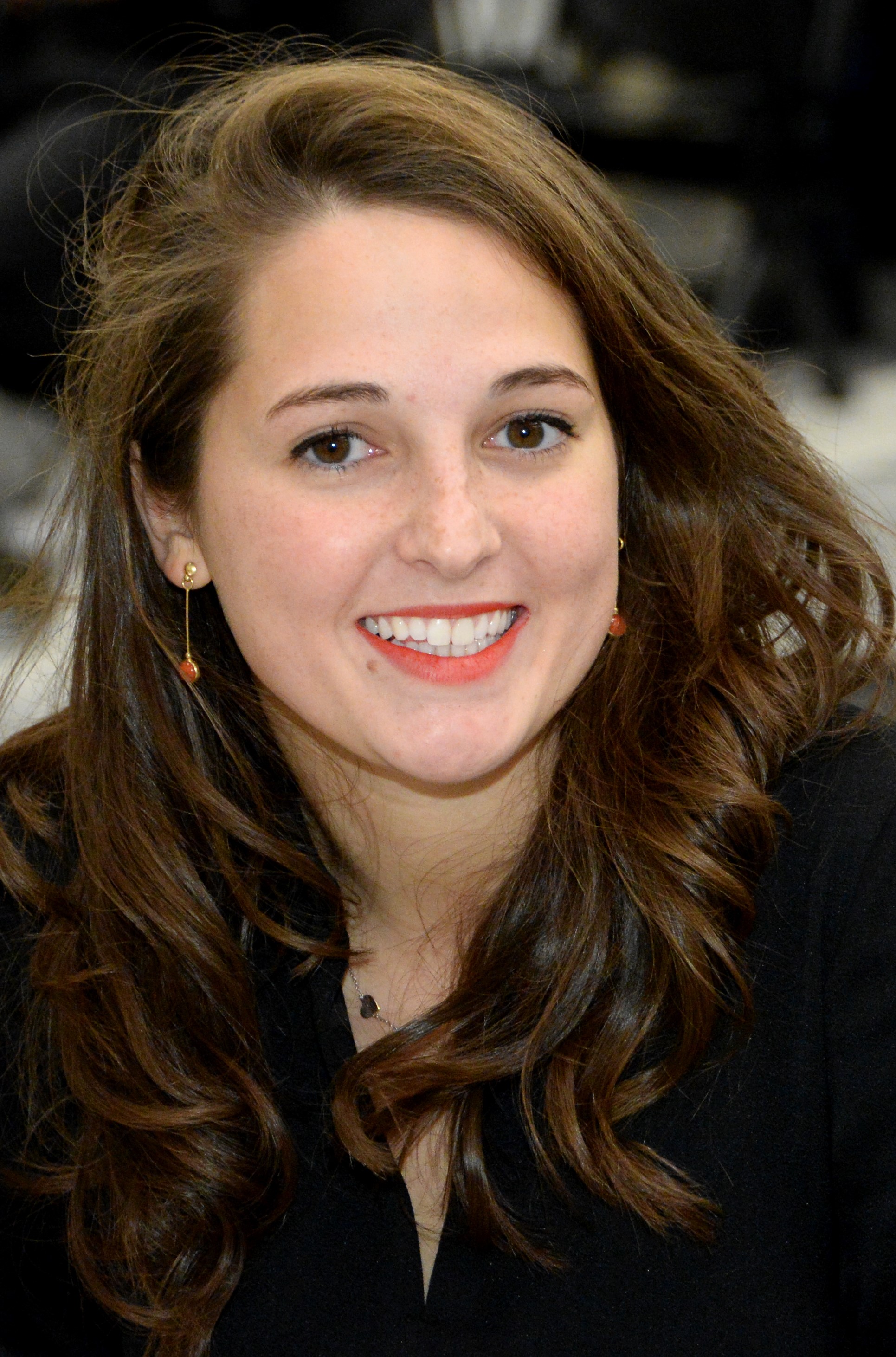}}]
	{Martina Vanelli} \tcb{is an Assistant Professor in the Faculty of Science and Engineering, University of Groningen. She} received the B.Sc., M.S., and Ph.D.~degrees in Applied Mathematics  from  Politecnico  di  Torino,  in  2017,  2019, and 2024, respectively. 
	\tcb{From 2024 to 2026, she was a Postdoctoral Fellow at the Institute for Information and Communication Technologies, Electronics and Applied Mathematics (ICTEAM), UCLouvain.} 
	Her research interests include identification, analysis, and control of multi-agent systems, with application to social and economic networks and power markets.
\end{IEEEbiography}

\begin{IEEEbiography}
	[{\includegraphics[width=1in,height=1.25in,clip,keepaspectratio]{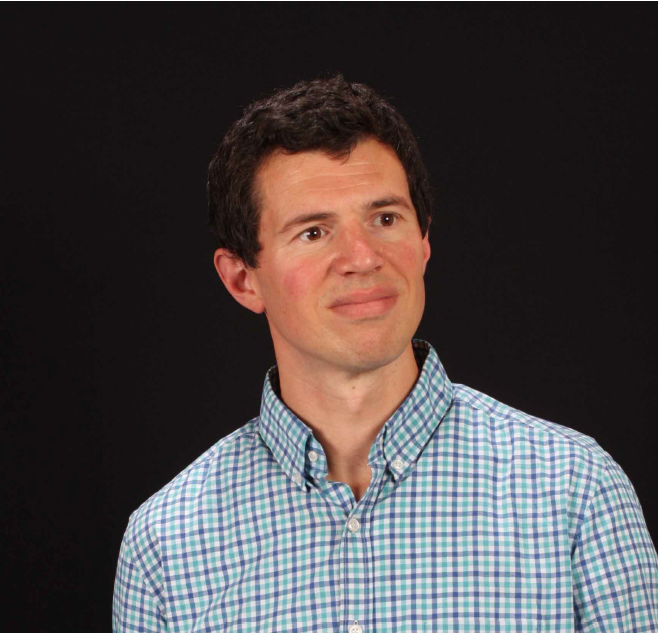}}]
	{Julien M. Hendrickx} is professor of mathematical
	engineering at UCLouvain, in the Ecole Polytechnique
	de Louvain since 2010. He obtained
	an engineering degree in applied mathematics
	(2004) and a PhD in mathematical engineering
	(2008) from the same university. He has been
	a visiting researcher at the University of Illinois
	at Urbana Champaign in 2003-2004, at the National
	ICT Australia in 2005 and 2006, and at the
	Massachusetts Institute of Technology in 2006
	and 2008. He was a postdoctoral fellow at the
	Laboratory for Information and Decision Systems of the Massachusetts	Institute of Technology 2009 and 2010, holding postdoctoral fellowships	of the F.R.S.- FNRS (Fund for Scientific Research) and of Belgian American Education Foundation.  He was also resident scholar at the Center	for Information and Systems Engineering (Boston University) in 2018-2019, holding a WBI.World excellence fellowship. Doctor Hendrickx is the recipient of the 2008 EECI award for the best PhD thesis in Europe	in the field of Embedded and Networked Control, and of the Alcatel-Lucent-Bell 2009 award for a PhD thesis on original new concepts or	application in the domain of information or communication technologies.
\end{IEEEbiography}
\vfill
		\end{document}